\input ./style/arxiv-general.cfg
\documentclass[aap,MSNbibl,seceqn,dvips]{arximspdf}
\makeatletter
   \@ifpackageloaded{graphicx}{}{\usepackage{graphicx}}
\makeatother

\doi{10.1214/14-AAP1054} 
\volume{25}
\issue{5}
\pubyear{2015}
\firstpage{2503}
\lastpage{2534}
\docsubty{FLA}

\makeatletter
\newcommand{\rrvert}{\vert}
\newcommand{\rrVert}{\Vert}
\newcommand{\llvert}{\vert}
\newcommand{\llVert}{\Vert}
\newtheorem{theorem}{Theorem}[section]
\newtheorem{proposition}[theorem]{Proposition}
\newtheorem{lemma}[theorem]{Lemma}
\newtheorem{corollary}[theorem]{Corollary}
\newtheorem{remark}[theorem]{Remark}
\newproclaim{example}[theorem]{Example}
\newproclaim{assumption}[theorem]{Assumption}
\newcommand{\inv}{\mathrm{inv}}
\newcommand{\mathbbm}{\mathbb}
\newcommand{\F}{\mathbb{F}}
\newcommand{\G}{\mathbb{G}}
\newcommand{\N}{\mathbb{N}}
\renewcommand{\P}{\mathbb{P}}
\newcommand{\Q}{\mathbb{Q}}
\newcommand{\R}{\mathbb{R}}
\renewcommand{\S}{\mathbb{S}}
\newcommand{\T}{\mathbb{T}}

\newcommand{\cA}{\mathcal{A}}
\newcommand{\cE}{\mathcal{E}}
\newcommand{\cF}{\mathcal{F}}
\newcommand{\cG}{\mathcal{G}}
\newcommand{\cH}{\mathcal{H}}
\newcommand{\cP}{\mathcal{P}}
\newcommand{\cT}{\mathcal{T}}
\newcommand{\cU}{\mathcal{U}}
\newcommand{\iint}{\int\!\!\int}
\newcommand{\1}{\mathbf{1}}
\newcommand{\bd}{\mathbf{d}}
\newcommand{\tomega}{\tilde{\omega}}
\newcommand{\bomega}{\bar{\omega}}
\newcommand{\btau}{{\bar{\tau}}}
\newcommand{\fPO}{\mathfrak{P}(\Omega)}
\makeatother

\begin{document}
\begin{frontmatter}

\title{Optimal stopping under adverse nonlinear expectation and
related games}
\runtitle{Optimal stopping under adverse nonlinear expectation}

\begin{aug}
\author[A]{\fnms{Marcel}~\snm{Nutz}\corref{}\ead[label=e1]{mnutz@math.columbia.edu}\thanksref{T1}}
\and
\author[B]{\fnms{Jianfeng}~\snm{Zhang}\ead[label=e2]{jianfenz@usc.edu}\thanksref{T2}}
\runauthor{M. Nutz and J. Zhang}
\affiliation{Columbia University and University of Southern California}
\address[A]{Department of Mathematics\\
Columbia University\\
2990 Broadway\\
New York, New York 10027\\
USA\\
\printead{e1}} 
\address[B]{Department of Mathematics\\
University of Southern California\\
3620 S. Vermont Avenue\\
Los Angeles, California 90089\\
USA\\
\printead{e2}}
\end{aug}
\thankstext{T1}{Supported by NSF Grant DMS-12-08985.}
\thankstext{T2}{Supported by NSF Grant DMS-10-08873.}

\received{\smonth{8} \syear{2013}}
\revised{\smonth{6} \syear{2014}}

%
\begin{abstract}
We study the existence of optimal actions in a zero-sum game $\inf_\tau\sup_P E^P[X_\tau]$ between a stopper and a controller choosing
a probability measure. This includes the optimal stopping problem $\inf_\tau\mathcal{E}(X_\tau)$ for a class of sublinear expectations
$\mathcal{E}(\cdot)$ such as the $G$-expectation. We show that the
game has a value. Moreover, exploiting the theory of sublinear
expectations, we define a nonlinear Snell envelope $Y$ and prove that
the first hitting time $\inf\{t\dvtx   Y_t=X_t\}$ is an optimal stopping
time. The existence of a saddle point is shown under
a compactness condition. Finally, the results are applied to the
subhedging of American options under volatility uncertainty.
\end{abstract}

%
\end{frontmatter}

\section{Introduction}\label{seintro}

On the space of continuous paths, we study a zero-sum stochastic game
%
\begin{equation}
\label{eqgameIntro} \inf_{\tau\in\cT} \sup_{P\in\cP}
E^P[X_\tau]
\end{equation}
between a stopper and a controller;
here, $X=(X_t)$ is the process to be stopped, $\cT$~is the set of
stopping times with values in a given interval $[0,T]$, and $\cP$ is a
given set of probability measures.
Specifically, we are interested in the situation where $\cP$ may be
nondominated, that is, there is no reference measure with respect to
which all $P\in\cP$ are absolutely continuous. This is the case, for
instance, when $\cP$ is the set of laws resulting from a controlled
stochastic differential equation whose diffusion coefficient is
affected by the control (cf. Example~\ref{excontrolledSDE}), whereas
the dominated case would correspond to the case where only the drift is
controlled. Or, in the language of partial differential equations, we
are interested in the fully nonlinear case rather than the semilinear
case. Technically, the nondominated situation entails that general
minimax results cannot be applied to~(\ref{eqgameIntro}), that the
cost functional $\sup_{P\in\cP} E^P[ \cdot ]$ does not satisfy
the dominated convergence theorem, and of course the absence of various
tools from stochastic analysis.
Our main results for the controller-and-stopper game include the
existence of an optimal action $\tau^*$ for the stopper under general
conditions and the existence of a saddle point $(\tau^*,P^*)$ under a
compactness condition. Both of these results were previously known only
in the case of drift control; cf. the review of literature at the end
of this section.

If we introduce the sublinear expectation $\cE(\cdot)=\sup_{P\in\cP
} E^P[ \cdot ]$, the stopper's part of the game can also be
interpreted as the nonlinear optimal stopping problem
%
\begin{equation}
\label{eqstopIntro} \inf_{\tau\in\cT} \cE(X_\tau).
\end{equation}
This alternate point of view is of independent interest, but it will
also prove to be useful in establishing the existence of optimal
actions for the game. Indeed, we shall start our analysis with~(\ref
{eqstopIntro}) and exploit the theory of nonlinear expectations, which
suggests to mimic the classical theory of optimal stopping under linear
expectation (e.g., \cite{ElKaroui81}). Namely, we define the Snell envelope
\[
Y_t= \inf_{\tau\in\cT_t} \cE_t(X_\tau),
\]
where $\cT_t$ is the set of stopping times with values in $[t,T]$ and
$\cE_t(\cdot)$ is the conditional sublinear expectation as obtained
by following the construction of~\cite{NutzVanHandel12}. Under
suitable assumptions, we show that the first hitting time
\[
\tau^*=\inf\{t\dvtx  Y_t=X_t\}
\]
is a stopping time which is optimal; that is, $\cE(X_{\tau^*})=\inf_{\tau\in\cT} \cE(X_\tau)$. Armed with this result, we return to
the game-theoretic point of view and
prove the existence of the value,
\[
\inf_{\tau\in\cT} \sup_{P\in\cP} E^P[X_\tau]
= \sup_{P\in\cP
} \inf_{\tau\in\cT}
E^P[X_\tau].
\]
Moreover, under a weak compactness assumption on $\cP$, we construct
$P^*\in\cP$ such that $(\tau^*,P^*)$ is a saddle point for~(\ref
{eqgameIntro}). These three main results are summarized in
Theorem~\ref{thoptimal}. Finally, we give an application to the
financial problem of pricing a path-dependent American option under
volatility uncertainty and show in Theorem~\ref{thduality} that~(\ref
{eqstopIntro}) yields the buyer's price (or subhedging price) in an
appropriate financial market model.

It is worth remarking that our results are obtained by working
``globally'' and not, as is often the case in the study of
continuous-time games, by a local-to-global passage based on a
Bellman--Isaacs operator; in fact, all ingredients of our setup can be
non-Markovian (i.e., path-dependent). The ``weak'' formulation of the
game, where the canonical process plays the role of the state process,
is important in this respect.

Like in the classical stopping theory, a dynamic programming principle
plays a key role in our analysis. We first prove this principle for the
upper value function $Y$ rather than the lower one (sup--inf), which
would be the standard choice in the literature on robust optimal
stopping (but of course the result for the lower value follows once the
existence of the value is established). The reason is that, due to the
absence of a reference measure, the structure of the set $\cT$ of
stopping times is inconvenient for measurable selections, whereas on
the set $\cP$ we can exploit the natural Polish structure. Once again
due to the absence of a reference measure, we are unable to infer the
optimality of $\tau^*$ directly from the dynamic programming
principle. However, we observe that in the discrete-time case, the
classical recursive analysis can be carried over rather easily by
exploiting the tower property of the nonlinear expectation. This
recursive structure
extends to the case where the processes are running in continuous time
but the stopping times are restricted to take values in a given
discrete set, like for a Bermudan option. To obtain the optimality of
$\tau^*$, we then approximate the continuous-time problem with such
discrete ones; the key idea is to compare the first hitting times for
the discrete problems with the times $\tau^\varepsilon=\inf\{t\dvtx
X_t-Y_t\leq\varepsilon\}$ and exploit the $\cE$-martingale property
of the discrete-time Snell envelope.
A similar approximation is also used to prove the existence of the
value, as it allows to circumvent the mentioned difficulty in working
with the lower value function: we first identify the upper and lower
value functions in the discrete problems and then pass to the limit.
Finally, for the existence of $P^*$, an important difficulty is that we
have little information about the regularity of $X_{\tau^*}$. Our
construction uses a compactness argument and a result of~\cite{EkrenTouziZhang12stop} on the approximation of hitting times by
random times that are continuous in $\omega$ to find a measure
$P^\sharp$ which is optimal up to the time $\tau^*$. In a second
step, we manipulate $P^\sharp$ in such a way that for the stopper,
immediate stopping after $\tau^*$ is optimal, which yields the optimal
measure $P^*$ for the full time interval. All this is quite different
from the existing arguments for the dominated case.

While the remainder of this \hyperref[seintro]{Introduction} concerns the related
literature, the rest of the paper is organized as follows: Section~\ref{seprelims} details the setup and the construction of the sublinear
expectation. In Section~\ref{semainResults}, we state our main
result, discuss its assumptions and give a concrete example related to
controlled stochastic functional/differential equations. Section~\ref{seproofs} contains the proof of the main result, while the
application to option pricing is studied in Section~\ref{seamericanOption}.

\subsection*{Literature}
In terms of the mathematics involved, the
study of the problem $\sup_{\tau\in\cT} \sup_{P\in\cP}
E^P[X_\tau]$ in \cite{EkrenTouziZhang12stop} is the closest to the
present one. Although this is a control problem with discretionary
stopping rather than a game, their regularity results are similar to
ours. On the other hand, the proofs of the optimality of $\tau^*$ are
completely different: in~\cite{EkrenTouziZhang12stop}, it was
relatively simple to obtain the martingale property up to $\tau
^\varepsilon$, directly in continuous time, and the main difficulty
was the passage from $\tau^\varepsilon$ to $\tau^*$, which is
trivial in our case. (The existence of an optimal $P^*\in\cP$ was not
studied in~\cite{EkrenTouziZhang12stop}.) Somewhat surprisingly, the
conditions obtained in the present paper are weaker than the ones
in~\cite{EkrenTouziZhang12stop}; in particular, for the optimality of
$\tau^*$, we do not assume that $\cP$ is compact.

After the publication of the preprint of the present work, \cite{BayraktarYao13} showed the existence of the optimal stopping time and
value in a case where $X$ is not bounded (under various other
assumptions). The authors go through the dynamic programming for the
lower value rather than the upper one, by using approximations based on
the regularity of $X$. The existence of a saddle point is not addressed
directly and we mention that our result does not apply, because the
technical assumptions of \cite{BayraktarYao13} preclude closedness of
$\cP$ in most cases of interest.

For the case where $\cP$ is dominated, the problem of optimal stopping
under nonlinear expectation (and related risk measures) is fairly well
studied; see, in particular, \cite{BayraktarKaratzasYao10,BayraktarYao11a,BayraktarYao11b,ElKarouiPardouxQuenez97,FollmerSchied04,KaratzasZamfirescu05,KaratzasZamfirescu08,Riedel09}.
The mathematical analysis for that case is quite different. On the
other hand, there is a literature on controller-and-stopper games. In
the discrete-time case, \cite{MaitraSudderth96} obtained a general
result on the existence of the value. For the continuous-time problem,
the literature is divided into two parts: in the non-Markovian case,
only the pure drift control has been studied, cf. \cite{KaratzasZamfirescu08} and the references therein; this again
corresponds to the dominated situation. For the nondominated situation,
results exist only in the Markovian case, where the presence of
singular measures plays a lesser role; cf. \cite{BayraktarHuang11,KaratzasSudderth01,KaratzasSudderth07,KaratzasZamfirescu08}.
In particular, \cite{BayraktarHuang11} obtained the existence of the
value for a diffusion setting via the comparison principle for the
associated partial differential equation. On the other hand, \cite{KaratzasSudderth01} studied a linear diffusion valued in the unit
interval with absorbing boundaries and found, based on scale-function
arguments, rather explicit formulas for the value and a saddle point.
Apart from such rather specific models, our results on the existence of
optimal actions are new even in the Markovian case.

Regarding the literature on nonlinear expectations, we refer to \cite{Peng07,Peng10icm} and the references therein; for the related
second-order backward stochastic differential equations (2BSDE) to
\cite{CheriditoSonerTouziVictoir07,SonerTouziZhang2010dual,SonerTouziZhang2010bsde}, and in particular to~\cite{MatoussiPiozinPossamai12,MatoussiPossamaiZhou12} for the reflected
2BSDE related to our problem; whereas for the uncertain volatility
model in finance, we refer to~\cite{AvellanedaLevyParas95,Lyons95,Smith02}.

\section{Preliminaries}\label{seprelims}

In this section, we introduce the setup and in particular the sublinear
expectation. We follow~\cite{NutzVanHandel12} as we need the
conditional expectation to be defined at every path and for all Borel functions.

\subsection{Notation}

We fix $d\in\N$ and let $\Omega=\{\omega\in C(\R_+;\R^d)\dvtx
\omega_0=0\}$ be the space of continuous paths equipped with the
topology of locally uniform convergence and the Borel $\sigma$-field
$\cF= \mathcal{B}(\Omega)$.
We denote by $B=(B_t)_{t \geq0}$ the canonical process $B_t(\omega
)=\omega_t$
and by $\F=(\cF_t)_{t \geq0}$
the (raw) filtration generated by $B$. Finally, $\mathfrak{P}(\Omega
)$ denotes the space of probability measures on $\Omega$ with the
topology of weak convergence.
Throughout this paper, ``stopping time'' will refer to a finite $\F
$-stopping time. Given a stopping time $\tau$ and $\omega, \omega
'\in\Omega$, we set
\[
\bigl(\omega\otimes_\tau\omega'\bigr)_u:=
\omega_u \1_{[0,\tau(\omega
))}(u) + \bigl(\omega_{\tau(\omega)} +
\omega'_{u-\tau(\omega
)} \bigr) \1_{[\tau(\omega), \infty)}(u),\qquad u\geq0.
\]
For any probability measure $P\in\mathfrak{P}(\Omega)$, there is a
regular conditional
probability distribution $\{P^\omega_\tau\}_{\omega\in\Omega}$
given $\cF_\tau$ satisfying
\[
P^\omega_\tau \bigl\{\omega'\in\Omega\dvtx
\omega' = \omega\mbox{ on } \bigl[0,\tau(\omega)\bigr] \bigr\} = 1
\qquad\mbox{for all }\omega\in \Omega;
\]
cf. \cite{StroockVaradhan79}, page~34. We then define $P^{\tau,\omega}\in\mathfrak{P}(\Omega)$ by
\[
P^{\tau,\omega}(A):=P^\omega_\tau(\omega\otimes_\tau
A),\qquad A\in\cF\mbox{ where }\omega\otimes_\tau A:=\bigl\{\omega
\otimes _\tau\omega'\dvtx  \omega'\in A\bigr\}.
\]
Given a function $f$ on $\Omega$ and $\omega\in\Omega$,
we also define the function $f^{\tau,\omega}$ by
\[
f^{\tau,\omega}\bigl(\omega'\bigr):=f\bigl(\omega
\otimes_\tau\omega'\bigr),\qquad \omega'\in
\Omega.
\]
We then have $E^{P^{\tau,\omega}}[f^{\tau,\omega}]=E^P[f| \cF_\tau
](\omega)$ for $P$-a.e. $\omega\in\Omega$.

\subsection{Sublinear expectation}

Let $\{\cP(s,\omega)\}_{(s,\omega)\in\R_+\times\Omega}$ be a
family of subsets of $\mathfrak{P}(\Omega)$, adapted in the sense that
\[
\cP(s,\omega)=\cP\bigl(s,\omega'\bigr)\qquad\mbox{if }\omega
|_{[0,s]}=\omega'| _{[0,s]},
\]
and define $\cP(\tau,\omega):=\cP(\tau(\omega),\omega)$ for any
stopping time $\tau$.
Note that the set $\cP(0,\omega)$ is independent of $\omega$ as all
paths start at the origin. Thus, we can define $\cP:=\cP(0,\omega)$.
We assume throughout that $\cP(s,\omega)\neq\varnothing$ for all
$(s,\omega)\in\R_+\times\Omega$.

The following assumption, which is in force throughout the paper, will
enable us to construct the conditional sublinear expectation related to
$\{\cP(s,\omega)\}$; it essentially states that our problem admits
dynamic programming. We recall that a subset of a Polish space is
called analytic if it is the image of a Borel subset of another Polish
space under a Borel mapping (we refer to \cite{BertsekasShreve78}, Chapter~7, for background).

\begin{assumption}\label{asanalytic}
Let $s\in\R_+$, let $\tau$ be a stopping time such that $\tau\geq
s$, let $\bomega\in\Omega$ and $P\in\cP(s,\bomega)$. Set $\theta:=\tau^{s,\bomega}-s$.
\begin{longlist}[(iii)]
\item The graph $\{(P',\omega)\dvtx  \omega\in\Omega,  P'\in\cP(\tau,\omega)\}  \subseteq  \fPO\times\Omega$ is analytic.

\item We have $P^{\theta,\omega} \in\cP(\tau,\bomega\otimes
_s\omega)$ for $P$-a.e. $\omega\in\Omega$.

\item If $\nu\dvtx  \Omega\to\fPO$ is an $\cF_\theta$-measurable
kernel and $\nu(\omega)\in\cP(\tau,\bomega\otimes_s\omega)$ for
$P$-a.e. $\omega\in\Omega$,
then the measure defined by
\[
\bar{P}(A)=\iint(\1_A)^{\theta,\omega}\bigl(\omega'
\bigr) \nu\bigl(d\omega ';\omega\bigr) P(d\omega),\qquad A\in\cF
\]
is an element of $\cP(s,\bomega)$.
\end{longlist}
\end{assumption}

Let us recall that a function $f\dvtx \Omega\to\overline{\R}$ is called
upper semianalytic if $\{f>c\}$ is analytic for each $c\in\R$; in
particular, every Borel function is upper semianalytic (cf. \cite{BertsekasShreve78}, Chapter~7).
Moreover, we recall that the universal completion of a $\sigma$-field
$\cA$ is given by $\cA^{*}:=\bigcap_{P} \cA^{P}$, where $\cA^{P}$
denotes the completion with respect to $P$ and the intersection is
taken over all probability measures on~$\cA$.
Let us agree that $E^P[f]:=-\infty$ if $E^P[f^+]=E^P[f^-]=+\infty$,
then we can introduce the sublinear expectation corresponding to $\{\cP
(s,\omega)\}$ as follows (cf. \cite{NutzVanHandel12}, Theorem~2.3).

\begin{proposition}\label{prNvH}
Let $\sigma\leq\tau$ be stopping times and let $f\dvtx \Omega\to
\overline{\R}$ be an upper semianalytic function. Then the function
%
\begin{equation}
\label{eqdefE} \cE_\tau(f) (\omega):=\sup_{P\in\cP(\tau,\omega)}
E^P\bigl[f^{\tau,\omega}\bigr],\qquad\omega\in\Omega
\end{equation}
is $\cF_\tau^*$-measurable and upper semianalytic. Moreover,
%
\begin{equation}
\label{eqtower} \cE_\sigma(f) (\omega) = \cE_\sigma\bigl(
\cE_\tau(f)\bigr) (\omega)\qquad\mbox{for all }\omega\in\Omega.
\end{equation}
\end{proposition}

We write $\cE(\cdot)$ for $\cE_0(\cdot)$.
We shall use very frequently (and often implicitly) the following
extension of Galmarino's test (cf. \cite{NutzVanHandel12}, Lemma~2.5).

\begin{lemma}\label{leunivgalmarino}
Let $f\dvtx \Omega\to\overline{\R}$ be $\cF^*$-measurable and let $\tau$
be a stopping time. Then $f$ is $\cF_\tau^*$-measurable if and only
if $f(\omega)=f(\omega_{\cdot\wedge\tau(\omega)})$ for all
$\omega\in\Omega$.
\end{lemma}

The following is an example for the use of Lemma~\ref{leunivgalmarino}: if $f$ and $g$ are bounded and upper semianalytic,
and $g$ is $\cF^*_t$ measurable, then
\[
\cE_t(f+g)=\cE_t(f)+g
\]
by~(\ref{eqdefE}), since the test shows that $g^{t,\omega}=g$.
Similarly, if we also have that $g\geq0$, then $\cE_t(fg)=\cE
_t(f)g$. We emphasize that all these identities hold at every single
$\omega\in\Omega$, without an exceptional set.

The most basic example we have in mind for $\cE(\cdot)$ is the
$G$-expectation of \mbox{\cite{Peng07,Peng08}}; or more precisely, its
extension to the upper semianalytic functions. In this case, $\cP
(s,\omega)$ is actually independent of $(s,\omega)$; more general
cases are discussed in Section~\ref{sesuffCond}.

\begin{example}[($G$-expectation)]\label{exGexp}
Let $U\neq\varnothing$ be a convex, compact set of nonnegative definite
symmetric $d\times d$ matrices and define $\cP_G$ to be the set of all
probabilities on $\Omega$ under which the canonical process $B$ is a
martingale whose quadratic variation $\langle B \rangle$ is absolutely
continuous $dt\times P$-a.e. and
\[
\frac{d\langle B \rangle_t}{dt}\in U\qquad dt\times P\mbox{-a.e.}
\]
Moreover, set $\cP(s,\omega):=\cP_G$ for all $(s,\omega)\in\R
_+\times\Omega$. Then Assumption~\ref{asanalytic} is satisfied
(cf. \cite{NutzVanHandel12}, Theorem~4.3) and $\cE(\cdot)$ is
called the $G$-expectation associated with $U$ (where $2G$ is the
support function of $U$). We remark that $\cP(s,\omega)$ is weakly
compact in this setup.

More generally, Assumption~\ref{asanalytic} is established in~\cite{NutzVanHandel12} when $U$ is a set-valued process (i.e., a Borel set
depending on $t$ and $\omega$). In this case, $\cP(s,\omega)$ need
not be compact and depends on $(s,\omega)$.
\end{example}

\section{Main results}\label{semainResults}

Let $T\in(0,\infty)$ be the time horizon. For any $t\in[0,T]$, we
denote by $\cT_t$ the set of all $[t,T]$-valued stopping times. For
technical reasons, we shall also consider the smaller set
$\cT^t\subseteq\cT_t$ of stopping times which do not depend on the
path up to time $t$; that is,
%
\begin{equation}
\label{eqdefTsupt} \cT^t=\bigl\{\tau\in\cT_t\dvtx
\tau^{t,\omega}=\tau^{t,\omega'}\mbox{ for all }\omega,
\omega'\in\Omega\bigr\}.
\end{equation}
In particular, $\cT:=\cT_0=\cT^0$ is the set of all $[0,T]$-valued
stopping times.
We define the pseudometric $\bd$ on $[0,T]\times\Omega$ by
\[
\bd\bigl[(t,\omega),\bigl(s,\omega'\bigr)\bigr]:=
\llvert t-s\rrvert + \bigl\llVert \omega_{\cdot\wedge
t}-\omega'_{\cdot\wedge s}
\bigr\rrVert _T
\]
%
for $(t,\omega), (s,\omega') \in[0,T]\times\Omega$, where $\llVert
\omega\rrVert  _u:=\sup_{r\leq u} \llvert  \omega_r\rrvert  $ for $u\geq0$ and $\llvert  \cdot\rrvert  $
is the Euclidean norm.

Let us now introduce the process $X$ to be stopped. Of course, the most
classical example is $X_t=f(B_t)$ for some function $f\dvtx \R^d\to\R$.
We consider a fairly general, possibly path-dependent functional
$X=X(B)$; note that the canonical process $B$ plays the role of the
state process. We shall work under the following regularity assumption.

\begin{assumption}\label{asXunifCont}
The process $X=(X_t)_{0\leq t\leq T}$ is progressively measurable,
uniformly bounded, has c\`adl\`ag trajectories with nonpositive jumps,
and there exists a modulus of continuity $\rho_X$ such that
%
\begin{equation}
\label{eqXunifCont} X_{t}(\omega)-X_{s}\bigl(
\omega'\bigr)\leq\rho_X \bigl(\bd\bigl[(t,\omega ),
\bigl(s,\omega'\bigr)\bigr] \bigr)\qquad\mbox{for all }s\leq t
\end{equation}
and $\omega,\omega'\in\Omega$.
\end{assumption}

The subsequent assumptions are stated in a form that is convenient for
the proofs; in that sense, they are in the most general form.
Sufficient conditions and examples will be discussed in Section~\ref{sesuffCond}.

\begin{assumption}\label{asEunifCont}
There is a modulus of continuity $\rho_{\cE}$ with the following
property. Let $t\in[0,T]$, $\tau\in\cT^t$ and $\bar\omega\in
\Omega$, then for all $\omega\in\Omega$ there exists $\tau_\omega
\in\cT^t$ such that
\[
\bigl\llvert \cE_t(X_\tau) (\bar\omega) -
\cE_t(X_{\tau_\omega}) (\omega)\bigr\rrvert \leq
\rho_{\cE}\bigl(\llVert \bar\omega-\omega\rrVert _t\bigr)
\]
and such that $(\omega,\omega')\mapsto\tau_\omega(\omega')$ is
$\cF_t\otimes\cF$-measurable.
\end{assumption}

We note that as $X$ is bounded, the moduli $\rho_X$ and $\rho_{\cE}$
can also be taken to be bounded.
In the subsequent assumption, we use the notation $B^\theta$ for the process
$B_{\cdot+\theta}-B_\theta$, where $\theta$ is a stopping time.

\begin{assumption}\label{asintegrability}
Let $\rho'$ be a bounded modulus of continuity. Then there exists a
modulus of continuity $\rho$ such that
\[
E^P \bigl[\rho'\bigl(\delta+ \bigl\llVert
B^\theta\bigr\rrVert _\delta\bigr) \bigr] \leq\rho(\delta ),
\qquad\delta\in[0,T] 
\]
for all $\theta\in\cT$, $P\in\cP(t,\omega)$ and $(t,\omega)\in
[0,T]\times\Omega$.
\end{assumption}

Let us now introduce the main object under consideration, the value
function (``nonlinear Snell envelope'') given by
%
\begin{equation}
\label{eqdefY} Y_t(\omega):= \inf_{\tau\in\cT_t}
\cE_t(X_\tau) (\omega),\qquad (t,\omega)\in[0,T]\times
\Omega.
\end{equation}
We shall see that $Y$ is Borel-measurable under the above assumptions,
and that $\cT_t$ can be replaced by $\cT^t$ without changing the
value of $Y_t$. The following is our main result.

\begin{theorem}\label{thoptimal}
Let Assumptions~\ref{asXunifCont},~\ref{asEunifCont} and~\ref{asintegrability} hold.
\begin{longlist}[(iii)]
\item There exists an optimal stopping time; namely,
\[
\tau^*:=\inf\bigl\{t\in[0,T]\dvtx  Y_t=X_t\bigr\}
\]
satisfies $\tau^*\in\cT$ and $\cE(X_{\tau^*})=\inf_{\tau\in\cT
} \cE(X_\tau)$.
\item The game has a value; that is,
\[
\inf_{\tau\in\cT}\sup_{P\in\cP} E^P[X_\tau]
= \sup_{P\in\cP} \inf_{\tau\in\cT}E^P[X_\tau].
\]
\item Suppose that $\cP(t,\omega)$ is weakly compact for all
$(t,\omega)\in[0,T]\times\Omega$. Then the game has a saddle point:
there exists $P^*\in\cP$ such that
\[
\inf_{\tau\in\cT}E^{P_*}[X_\tau] =
E^{P^*}[X_{\tau^*}] = \sup_{P\in\cP}
E^P[X_{\tau^*}].
\]
\end{longlist}
\end{theorem}

Of course, weak compactness in~(iii) refers to the topology induced by
the continuous bounded functions, and a similar identity as in (ii)
holds for the value functions at positive times (cf. Lemma~\ref
{leexistenceOfValue}). The proof of the theorem is stated in
Section~\ref{seproofs}. We mention the following variant of
Theorem~\ref{thoptimal}(i) where the stopping times are restricted to
take values in a discrete set $\T$; in this case, Assumption~\ref
{asintegrability} is unnecessary. The proof is again deferred to the
subsequent section.

\begin{remark}\label{rkoptimalDiscrete}
Let Assumptions~\ref{asXunifCont} and~\ref{asEunifCont} hold. Let
$\T=\{t_0,t_1,\ldots,t_n\}$, where $n\in\N$ and $t_0<t_1<\cdots
<t_n=T$, and consider the obvious corresponding notions like
$\cT_t(\T)=\{\tau\in\cT_t\dvtx \tau(\cdot)\in\T\}$ and $Y_t=\inf_{\tau\in\cT_t(\T)} \cE_t(X_\tau)$.
Then $Y$ satisfies the backward recursion
\[
Y_{t_n}=X_{t_n}\quad\mbox{and}\quad Y_{t_i} =
X_{t_i} \wedge\cE _{t_i}(Y_{t_{i+1}}),\qquad i=0,\ldots,n-1
\]
and $\tau^*:=\inf\{t\in\T\dvtx   Y_t=X_t\}$ satisfies $\cE(X_{\tau
^*})=\inf_{\tau\in\cT(\T)} \cE(X_\tau)$.
\end{remark}

\subsection{Sufficient conditions for the main assumptions}\label{sesuffCond}

In the remainder of this section, we discuss the conditions of the
theorem. Assumption~\ref{asXunifCont} is clearly satisfied when
$X_t=f(B_t)$ for a bounded, uniformly continuous function $f$. The
following shows that Assumption~\ref{asEunifCont} is trivially
satisfied, for example, for the \mbox{$G$-}expectation of Example~\ref
{exGexp} (a nontrivial situation is discussed in Example~\ref
{excontrolledSDE}).

\begin{remark}\label{rkEunifContForIndep}
Assume that $\cP(t,\omega)$ does not depend on $\omega$, for all
$t\in[0,T]$. Then Assumption~\ref{asXunifCont} implies
Assumption~\ref{asEunifCont}.
\end{remark}

\begin{pf}
Let $\tau\in\cT^t$, then $\tau^{t,\omega}=\tau^{t,\omega
'}=:\theta$ for all $\omega,\omega'\in\Omega$. Moreover, taking
$s=t$ in~(\ref{eqXunifCont}) shows that
\[
\bigl\llvert X_{s}(\omega)-X_{s}\bigl(
\omega'\bigr)\bigr\rrvert \leq\rho_X\bigl(\bigl\llVert
\omega-\omega'\bigr\rrVert _s\bigr).
\]
We deduce that for all $\tomega\in\Omega$,
%
\begin{eqnarray} \label{eqproofEunifContForIndep}
\bigl\llvert (X_\tau)^{t,\omega}(\tomega)-(X_\tau)^{t,\omega'}(
\tomega)\bigr\rrvert &=& \bigl\llvert X_{\theta(\tomega)}(\omega
\otimes_t\tomega) - X_{\theta
(\tomega)}\bigl(\omega'
\otimes_t\tomega\bigr)\bigr\rrvert
\nonumber
\\
& \leq&\rho_X\bigl(\bigl\llVert \omega\otimes_t \tomega-
\omega'\otimes_t \tomega \bigr\rrVert _{\theta(\tomega)}
\bigr)
\\
&=& \rho_X\bigl(\bigl\llVert \omega-\omega'\bigr
\rrVert _t\bigr).\nonumber
\end{eqnarray}
If $\cP(t,\cdot)=\cP(t)$, it follows that
\begin{eqnarray*}
\bigl\llvert \cE_t(X_\tau) (\omega) -
\cE_t(X_\tau) \bigl(\omega'\bigr)\bigr\rrvert
&\leq& \sup_{P\in\cP(t)} E^P\bigl[\bigl\llvert
(X_\tau)^{t,\omega}-(X_\tau )^{t,\omega'}\bigr\rrvert
\bigr]
\\
&\leq&\rho_X\bigl(\bigl\llVert \omega-\omega'\bigr
\rrVert _t\bigr);
\end{eqnarray*}
that is, Assumption~\ref{asEunifCont} holds with $\rho_{\cE}=\rho
_X$ and $\tau_{\omega}=\tau$.
\end{pf}

The following is a fairly general sufficient condition for
Assumption~\ref{asintegrability}; it covers most cases of interest.

\begin{remark}\label{rkintegrability}
Suppose that for some $\alpha,c>0$, the moment condition
\[
E^P\bigl[\bigl\llVert B^\theta\bigr\rrVert _{\delta}
\bigr] \leq c \delta^\alpha,\qquad\delta\in[0,T]
\]
is satisfied for all $\theta\in\cT$, $P\in\cP(t,\omega)$ and
$(t,\omega)\in[0,T]\times\Omega$. Then Assumption~\ref
{asintegrability} holds.
In particular, this is the case if every $P\in\cP(t,\omega)$ is the
law of an It\^o process
\[
\Gamma_s= \int_0^s
\mu_r \,dr + \int_0^s
\sigma_r \,dW_r
\]
(where $W$ is a Brownian motion) and $\llvert  \mu\rrvert   + \llvert  \sigma\rrvert  \leq C$ for a
universal constant~$C$.
\end{remark}

\begin{pf}
Let $r\in\R$ be such that $\rho'\leq r$. For any $a>0$, we have
\[
E^P \bigl[\rho'\bigl(\delta+ \bigl\llVert
B^\theta\bigr\rrVert _{\delta}\bigr) \bigr] \leq\rho
'(\delta+ a) + r P\bigl\{\bigl\llVert B^\theta\bigr\rrVert
_{\delta}\geq a\bigr\}.
\]
Using that $P\{\llVert  B^\theta\rrVert  _{\delta}\geq a\} \leq a^{-1} E^P[\llVert
B^\theta\rrVert  _{\delta}] \leq a^{-1}c\delta^\alpha$ and choosing
$a=\delta^{\alpha/2}$, we obtain that
\[
E^P \bigl[\rho'\bigl(\delta+ \bigl\llVert
B^\theta\bigr\rrVert _{\delta}\bigr) \bigr] \leq\rho
'\bigl(\delta+ \delta^{\alpha/2}\bigr) + cr
\delta^{\alpha/2} =:\rho(\delta),
\]
which was the first claim.

Suppose that $B=A+M$ under $P$, where $\llvert  dA\rrvert  +d\llvert  \langle M \rangle\rrvert  \leq
C \,dt$; we focus on the scalar case for simplicity. Using the
Burkholder--Davis--Gundy inequalities for $M^\theta=M_{\cdot+\theta
}-M_\theta$, we have
\begin{eqnarray*}
E^P\bigl[\bigl\llVert B^\theta\bigr\rrVert _\delta
\bigr] &\leq& E^P\bigl[\bigl\llVert A^\theta\bigr\rrVert
_\delta\bigr]+E^P\bigl[\bigl\llVert M^\theta\bigr
\rrVert _\delta\bigr]
\\
&\leq& C\delta+c_1E^P\bigl[\bigl\llVert \bigl\langle
M^\theta \bigr\rangle^{1/2}\bigr\rrVert _\delta\bigr]
\\
&\leq&\bigl(CT^{1/2} + c_1C^{1/2}\bigr)
\delta^{1/2},\qquad\delta\in[0,T],
\end{eqnarray*}
where $c_1>0$ is a universal constant.
\end{pf}

Let us now discuss two classes of models. The first one is the main
example for control in the ``weak formulation,'' that is, the set of
controls is stated directly in terms of laws.

\begin{example}\label{exLevy}
Let $U$ be a nonempty, bounded Borel set of $\R^d\times\S_+$, where~$\S_{+}$ is the set of $d\times d$ nonnegative definite symmetric
matrices. Moreover, let $\cP$ be the set of all laws of continuous
semimartingales whose characteristics are absolutely continuous (with
respect to the Lebesgue measure) and whose \mbox{differential} characteristics
take values in $U$. That is, $\cP$ consists of laws of It\^o processes
$\int b_t \,dt + \int\sigma_t  \,dW_t$, each situated on its own
probability space with a Brownian motion $W$, where the pair $(b,\sigma
\sigma^\top)$ almost surely takes values in the set $U$. For
instance, if $d=1$ and $U=I_1\times I_2$ is a product of intervals,
this models the case where the controller can choose the drift from
$I_1$ and the (squared) diffusion from~$I_2$.

The above setup (and its extension to jump processes) is studied
in~\cite{NeufeldNutz13b}, where it is shown in particular that
Assumption~\ref{asanalytic} holds. Moreover, we see from Remark~\ref
{rkEunifContForIndep} that Assumption~\ref{asEunifCont} is
satisfied, while Remark~\ref{rkintegrability} shows that
Assumption~\ref{asintegrability} holds as well. Finally, if $U$ is
compact and convex, standard results (see~\cite{Zheng85}) imply that
the set $\cP$ is weakly compact.
\end{example}

In the second class of models, whose formulation is borrowed from~\cite{Nutz11}, the elements of $\cP$ correspond to the possible laws of
the solution to a controlled stochastic functional/differential
equation (SDE). This is the main case of interest for the
controller-and-stopper games in the ``strong formulation'' of control
and as we shall see, the sets $\cP(t,\omega)$ indeed depend on
$(t,\omega)$. Note that in the setting of the strong formulation the
set $\cP$ is typically not closed and in particular not compact, so
that we cannot expect the existence of a saddle point in general. For
simplicity, we only discuss the case of a driftless SDE.

\begin{example}\label{excontrolledSDE}
Let $U$ be a nonempty Borel set of $\R^d$ and let $\cU$ be the set of
all $U$-valued, progressively measurable, c\`adl\`ag processes $\nu$.
We denote by $\S_{++}$ the set of positive definite symmetric matrices
and by $\mathbbm{D}$ the Skorohod space of c\`adl\`ag paths in $\R
^{d}$ starting at the origin, and consider a function
\[
\sigma\dvtx  \R_+\times\mathbbm{D}\times U\to\S_{++}
\]
such that $(t,\omega)\mapsto\sigma(t,\Gamma(\omega),\nu_t(\omega
))$ is progressively measurable (c\`adl\`ag)
whenever $\Gamma$ and $\nu$ are progressively measurable (c\`adl\`ag).
We assume that $\sigma$ is uniformly Lipschitz in its second variable
with respect to the supremum norm,
and (for simplicity) uniformly bounded. Moreover, we assume that
$\sigma$ is a one-to-one function in its third variable, admitting a
measurable inverse on its range. More precisely, there exists a
function $\sigma^{\inv}\dvtx  \R_+\times\mathbbm{D}\times\S_{++}\to U$
such that
\[
\sigma^{\inv}\bigl(t,\omega, \sigma(t,\omega,u)\bigr) =u
\]
for all $(t,\omega,u)\in \R_+\times\mathbbm{D}\times U$, and
$\sigma^{\inv}$ satisfies the same measurability and c\`adl\`ag
properties as $\sigma$. Given $\nu\in\cU$, we consider the
stochastic functional/differential equation
\[
\Gamma_t=\int_0^t
\sigma(r,\Gamma,\nu_r) \,dB_r,\qquad t\geq0
\]
under the Wiener measure $P_0$ (i.e., $B$ is a $d$-dimensional Brownian
motion). This equation has a $P_0$-a.s. unique strong solution whose
law is denoted by $P(\nu)$. We then define $\cP=\{P(\nu)\dvtx   \nu\in
\cU\}$.
More generally, $\cP(s,\omega)$ is defined as the set of laws
$P(s,\omega,\nu)$ corresponding to the SDE with conditioned
coefficient $(r,\omega',u)\mapsto\sigma(r+s,\omega\otimes_s \omega
', u)$ and initial condition $\Gamma_0=\omega_s$; more precisely,
$P(s,\omega,\nu)$ is the law on $\Omega$ of the solution translated
to start at the origin (see also~\cite{Nutz11}).

In this model, Assumption~\ref{asanalytic} can be verified by the
arguments used in~\cite{Nutz11} and~\cite{NeufeldNutz12}; the
details are lengthy but routine. Assumption~\ref{asintegrability} is
satisfied by Remark~\ref{rkintegrability} since $\sigma$ is bounded
(note that this condition can be improved by using SDE estimates). We
impose Assumption~\ref{asXunifCont} on $X$ and turn to
Assumption~\ref{asEunifCont}, which is the main interest in this
example. The main problem in this regard is that we cannot impose
continuity conditions on the stopping times.

\begin{lemma}\label{lemasssss}
Assumption~\ref{asEunifCont} is satisfied in the present setting.
\end{lemma}

While we defer the actual proof to the \hyperref[app]{Appendix}, we sketch here the
rough idea for the case where there is no control in the SDE (i.e., $U$
is a singleton) and thus each set $\cP(s,\omega)$ consists of a
single measure $P(s,\omega)$.
Let $t\in[0,T]$, $\tau\in\cT^t$ and $\bar\omega\in\Omega$; we
shall construct $\tau_\omega\in\cT^t$ such that
\[
\bigl\llvert \cE_t(X_\tau) (\bar\omega) -
\cE_t(X_{\tau_\omega}) (\omega)\bigr\rrvert \leq
\rho_{\cE}\bigl(\llVert \bar\omega-\omega\rrVert _t\bigr)
\]
for all $\omega\in\Omega$ [and such that $(\omega,\omega')\mapsto
\tau_\omega(\omega')$ is $\cF_t\otimes\cF$-measurable].

Fix $\omega,\bar\omega\in\Omega$ and denote by $\Gamma^{t,\omega
}$ and $\Gamma^{t,\bar\omega}$ the corresponding solutions of the
SDE (translated to start at the origin), that is, we have
\begin{eqnarray}
d\Gamma^{t,\omega}_u = \sigma\bigl(u+t,\omega
\otimes_t \Gamma ^{t,\omega}\bigr) \,dB_u\quad\mbox{and}\quad d\Gamma^{t,\bar\omega}_u = \sigma\bigl(u+t,\bar\omega
\otimes_t \Gamma^{t,\bar\omega}\bigr) \,dB_u\nonumber
\\
\eqntext{P_0\mbox{-a.s.}}
\end{eqnarray}
Our aim is to construct $\tau_\omega\in\cT^t$ with the property that
%
\begin{equation}
\label{eqstoppingTransformSketch} \tau_\omega\bigl(0\otimes_t
\Gamma^{t,\omega}\bigr) = \tau\bigl(0\otimes _t
\Gamma^{t,\bar\omega}\bigr)\qquad P_0\mbox{-a.s.},
\end{equation}
where $0$ is the constant path (or any other path, for that matter).
Indeed, if this identity holds, then
\begin{eqnarray*}
\cE_t(X_{\tau_\omega}) (\omega) &=& E^{P(t,\omega)} \bigl[
(X_{\tau_\omega})^{t,\omega} \bigr]
\\
&=& E^{P(t,\omega)} \bigl[ X_{\tau_\omega(0\otimes_t \cdot)} (\omega \otimes_t
\cdot) \bigr]
\\
&=& E^{P_0} \bigl[ X_{\tau_\omega(0\otimes_t \Gamma^{t,\omega})} \bigl(\omega\otimes_t
\Gamma^{t,\omega}\bigr) \bigr]
\\
&=& E^{P_0} \bigl[ X_{\tau(0\otimes_t \Gamma^{t,\bar\omega})} \bigl(\omega\otimes_t
\Gamma^{t,\omega}\bigr) \bigr]
\end{eqnarray*}
and thus
%
\begin{eqnarray} \label{eqfinalEstimateSketch}
\qquad &&\bigl\llvert \cE_t(X_\tau) (\bar\omega) -
\cE_t(X_{\tau_\omega}) (\omega )\bigr\rrvert
\nonumber
\\
&&\qquad = \bigl\llvert E^{P_0} \bigl[ X_{\tau(0\otimes_t \Gamma^{t,\bar\omega})} \bigl(\bar\omega
\otimes_t \Gamma^{t,\bar\omega}\bigr) \bigr] - E^{P_0} \bigl[
X_{\tau(0\otimes_t \Gamma^{t,\bar\omega})} \bigl(\omega\otimes_t \Gamma^{t,\omega}\bigr)
\bigr] \bigr\rrvert
\nonumber
\\
&&\qquad \leq E^{P_0} \bigl[ \bigl\llvert X_{\tau(0\otimes_t \Gamma^{t,\bar\omega
})} \bigl(\bar\omega
\otimes_t \Gamma^{t,\bar\omega}\bigr) - X_{\tau
(0\otimes_t \Gamma^{t,\bar\omega})} \bigl(\omega
\otimes_t \Gamma ^{t,\omega}\bigr) \bigr\rrvert \bigr]
\\
&&\qquad \leq E^{P_0} \bigl[ \rho_{X} \bigl(\bigl\llVert \bar
\omega\otimes_t \Gamma^{t,\bar\omega} - \omega\otimes_t
\Gamma^{t,\omega}\bigr\rrVert _{T} \bigr) \bigr]
\nonumber
\\
&&\qquad \leq\rho_X\bigl(C \llVert \bar\omega-\omega \rrVert
_t\bigr),\nonumber
\end{eqnarray}
where the last inequality follows by a standard SDE estimate as
in~\cite{Nutz11}, Lemma~2.6, with $C>0$ depending only on the
Lipschitz constant of $\sigma$ and the time horizon $T$. This is the
desired estimate with $\rho_{\cE}(\cdot)=\rho_{X}(C\cdot)$.

To construct $\tau_\omega$ satisfying~(\ref
{eqstoppingTransformSketch}), we basically require a transformation
$\zeta^{\omega}\dvtx  \Omega\to\Omega$ mapping the paths of $\Gamma
^{t,\omega}$ to the corresponding paths of $\Gamma^{t,\bar\omega}$.
(The dependence of $\zeta^{\omega}$ on the fixed path $\bar\omega$
is suppressed in our notation.) Roughly speaking, this is accomplished
by the solution of
\[
\zeta= \int_0^{\cdot} \sigma(u+t,\bar\omega
\otimes_t \zeta )\sigma(u+t,\omega\otimes_t
B)^{-1} \,dB_u.
\]
Indeed, let us suppose for the moment that a solution $\zeta^{\omega
}$ can be defined in some meaningful way and\vspace*{1pt} that all paths of $\zeta
^{\omega}$ are continuous.
Then, formally, we have $\zeta^{\omega}(\Gamma^{t,\omega})=\Gamma
^{t,\bar\omega}$ and
%
\begin{equation}
\label{eqdefStoppingTransformSketch} \tau_\omega\bigl(\omega'\bigr):= \tau
\bigl(0\otimes_t \zeta^{\omega}\bigl(\omega '_{\cdot+t}-
\omega'_t\bigr)\bigr)
\end{equation}
defines a stopping time with the desired property~(\ref
{eqstoppingTransformSketch}). In the \hyperref[app]{Appendix}, we show how to make
this sketch rigorous and include the case of a controlled equation.
\end{example}

\section{Proof of Theorem~\texorpdfstring{\protect\ref{thoptimal}}{3.4}}\label{seproofs}

Assumptions~\ref{asXunifCont},~\ref{asEunifCont} and~\ref{asintegrability} are in force throughout this section, in which we
state the proof of Theorem~\ref{thoptimal} through a sequence of lemmas.

\subsection{Optimality of \texorpdfstring{$\tau^*$}{tau*}}

We begin with the optimality of $\tau^*$. All results have their
obvious analogues for the discrete case discussed in Remark~\ref
{rkoptimalDiscrete}; we shall state this separately only where necessary.
We first show that $\cT_t$ may be replaced by the set $\cT^t$
from~(\ref{eqdefTsupt}) in the definition~(\ref{eqdefY}) of $Y_t$.

\begin{lemma}\label{leYaltDef}
Let $t\in[0,T]$. Then
%
\begin{equation}
\label{eqYaltDef} Y_t= \inf_{\tau\in\cT^t}
\cE_t(X_\tau).
\end{equation}
\end{lemma}

\begin{pf}
The inequality ``$\leq$'' follows from the fact that $\cT^t\subseteq
\cT_t$. To see the reverse inequality, fix $\omega\in\Omega$ and
let $\varepsilon>0$. By the definition~(\ref{eqdefY}) of $Y_t$,
there exists $\tau\in\cT_t$ (depending on $\omega$) such that
\[
Y_t(\omega) \geq\cE_t(X_\tau) (\omega) -
\varepsilon.
\]
Define $\theta=\tau^{t,\omega}(B^t)$, where $B^t=B_{\cdot+t}-B_t$.
Clearly, $\theta\geq t$, and we see from Galmarino's test that $\theta
$ is a stopping time. Moreover, as a function of $B^t$, $\theta$ is
independent of the path up to time $t$; that is, $\theta\in\cT^t$.
Noting that $\tau^{t,\omega}=\theta^{t,\omega}$, the
definition~(\ref{eqdefE}) of $\cE_t(\cdot)$ shows that
$\cE_t(X_\tau)(\omega)=\cE_t(X_\theta)(\omega)$, and hence
\[
Y_t(\omega) \geq\cE_t(X_\theta) (\omega) -
\varepsilon.
\]
The result follows as $\varepsilon>0$ was arbitrary.
\end{pf}

\begin{lemma}\label{leYunifcont}
We have
%
\begin{equation}
\label{eqYunifCont} \bigl\llvert Y_t(\omega)-Y_{t}\bigl(
\omega'\bigr)\bigr\rrvert \leq\rho_{\cE}\bigl(\bigl\llVert
\omega-\omega'\bigr\rrVert _t\bigr)
\end{equation}
for all $t\in[0,T]$ and $\omega,\omega'\in\Omega$. In particular,
$Y_t$ is $\cF_t$-measurable.
\end{lemma}

\begin{pf}
In view of~(\ref{eqYaltDef}), this follows from Assumption~\ref{asEunifCont}.
\end{pf}

The following dynamic programming principle is at the heart of this section.

\begin{lemma}\label{leDPPdetTimes}
Let $0\leq s\leq t \leq T$. Then
\[
Y_s = \inf_{\tau\in\cT^s} \cE_s (
X_\tau\1_{\{\tau<t\}} + Y_t \1_{\{\tau\geq t\}} );
\]
moreover, $\cT^s$ can be replaced by $\cT_s$.
\end{lemma}

\begin{pf}
We first show the inequality ``$\geq$.'' Let $\tau\in\cT_s$. As
$\tau\vee t\in\cT_t$, we have $\cE_t(X_{\tau\vee t})\geq Y_t$ by
the definition~(\ref{eqdefY}) of $Y$. Using the tower property~(\ref
{eqtower}), it follows that
\[
\cE_s(X_\tau) = \cE_s \bigl( X_\tau
\1_{\{\tau<t\}} + \cE_t(X_{\tau
\vee t}) \1_{\{\tau\geq t\}}
\bigr) \geq\cE_s ( X_\tau\1_{\{
\tau<t\}} +
Y_t \1_{\{\tau\geq t\}} );
\]
where we have used Lemma~\ref{leunivgalmarino} and that all the
involved random variables are upper semianalytic. In view of~(\ref
{eqYaltDef}), taking the
infimum over $\tau\in\cT_s$ (resp., $\cT^s$) yields the claimed inequality.

We now turn to the inequality ``$\leq$.''
Fix $\tau\in\cT_s$ and set $\Lambda_0:=\{\tau<t\}$. Moreover, let
$\varepsilon>0$ and let $(\Lambda^i)_{i\geq1}$ be an $\cF
_t$-measurable partition of the set $\{\tau\geq t\}\in\cF_t$
such that the $\llVert  \cdot\rrVert  _t$-diameter of $\Lambda^i$ is smaller than
$\varepsilon$ for all $i\geq1$. Fix $\omega^i\in\Lambda^i$.
By~(\ref{eqYaltDef}), there exist $\tau^i\in\cT^t$ such that
\[
Y_t\bigl(\omega^i\bigr)\geq\cE_t(X_{\tau^i})
\bigl(\omega^i\bigr)-\varepsilon,\qquad i\geq1.
\]
In view of Assumption~\ref{asEunifCont} and~(\ref{eqYunifCont}),
there exist stopping times $\tau^i_\omega\in\cT^t$ such that
%
\begin{equation}
\label{eqproofDppSemicont} Y_t(\omega)\geq\cE_t(X_{\tau^i_\omega}) (
\omega) - \rho (\varepsilon),\qquad\omega\in\Lambda^i, i\geq1,
\end{equation}
where $\rho(\varepsilon)=\varepsilon+2\rho_{\cE}(\varepsilon)$.
Define $\hat{\tau}^i(\omega):= \tau^i_{\omega}(\omega)$ and
\[
\btau:= \tau\1_{\{\tau<t\}} + \sum_{i\geq1} \hat{
\tau}^i\1 _{\Lambda^i}.
\]
In view of the measurability condition in Assumption~\ref{asEunifCont}, we then have $\btau\in\cT_s$ and $\btau^{t,\omega
}=(\hat{\tau}^i)^{t,\omega} = (\tau^i_\omega)^{t,\omega}$ for
$\omega\in\Lambda^i$. Using also~(\ref{eqproofDppSemicont}), the
tower property, and $\{\tau<t\}=\{\btau<t\}$, we deduce that
\begin{eqnarray*}
\cE_s ( X_\tau\1_{\{\tau<t\}} + Y_t
\1_{\{\tau\geq t\}} ) &\geq&\cE_s \biggl( X_\tau
\1_{\{\tau<t\}} + \sum_{i\geq1} \cE
_t(X_{\hat{\tau}^i})\1_{\Lambda^i} \biggr) - \rho(\varepsilon)
\\
&=& \cE_s \biggl( X_\tau\1_{\{\tau<t\}} + \sum
_{i\geq1} \cE _t(X_{\btau})\1_{\Lambda^i}
\biggr) - \rho(\varepsilon)
\\
&=& \cE_s \bigl( X_\tau\1_{\{\tau<t\}} +
\cE_t(X_{\btau})\1_{\{\tau
\geq t\}} \bigr) - \rho(\varepsilon)
\\
&=& \cE_s ( X_\tau\1_{\{\tau<t\}} + X_{\btau}
\1_{\{\tau\geq t\}
} ) - \rho(\varepsilon)
\\
&=& \cE_s( X_\btau) - \rho(\varepsilon)
\\
&\geq& Y_s - \rho(\varepsilon).
\end{eqnarray*}
As $\tau\in\cT_s\supseteq\cT^s$ was arbitrary, the result follows
by letting $\varepsilon$ tend to zero.
\end{pf}

Based on the dynamic programming principle of Lemma~\ref
{leDPPdetTimes} and Assumption~\ref{asintegrability}, we can now
establish the path regularity of $Y$. The following is quite similar to
\cite{EkrenTouziZhang12stop}, Lemma~4.2.

\begin{lemma}\label{leYContTime}
There exists a modulus of continuity $\rho_Y$ such that
\[
\bigl\llvert Y_s(\omega)-Y_t(\omega)\bigr
\rrvert \leq\rho_Y\bigl(\bd\bigl[(s,\omega),(t,\omega)\bigr]\bigr)
\]
%
for all $s, t \in[0,T]$ and $\omega\in\Omega$.
\end{lemma}

\begin{pf}
We may assume that $s\leq t$. Using Lemma~\ref{leDPPdetTimes}, we have
\[
Y_s(\omega)-Y_t(\omega) \leq\cE_s(Y_t)
(\omega)-Y_t(\omega) \leq\cE_s \bigl( \bigl\llvert
Y_t -Y_t(\omega)\bigr\rrvert \bigr) (\omega),\qquad\omega
\in \Omega.
\]
[The right-hand side is the conditional expectation of the random
variable $\omega'\mapsto\llvert  Y_t(\omega') -Y_t(\omega)\rrvert  $, evaluated at
the point $\omega$.] On the other hand, Lemma~\ref{leDPPdetTimes}
and the subadditivity of $\cE_s(\cdot)$ yield that
\begin{eqnarray*}
Y_t(\omega)-Y_s(\omega) &=& Y_t(\omega) -
\inf_{\tau\in\cT^s} \cE_s ( X_\tau
\1_{\{
\tau<t\}} + Y_t \1_{\{\tau\geq t\}} ) (\omega)
\\
&\leq&\sup_{\tau\in\cT^s} \cE_s \bigl( Y_t(
\omega) -X_\tau\1_{\{
\tau<t\}} - Y_t \1_{\{\tau\geq t\}}
\bigr) (\omega)
\\
&=& \sup_{\tau\in\cT^s} \cE_s \bigl( Y_t(
\omega) -Y_t + (Y_t -X_\tau) \1_{\{\tau<t\}}
\bigr) (\omega)
\\
&\leq&\cE_s \bigl( \bigl\llvert Y_t -Y_t(
\omega)\bigr\rrvert \bigr) (\omega) + \sup_{\tau
\in\cT^s}
\cE_s \bigl((Y_t -X_\tau) \1_{\{\tau<t\}}
\bigr) (\omega).
\end{eqnarray*}
Combining these two estimates, we obtain that
%
\begin{equation}\label{eqproofYContTime}
\qquad\bigl\llvert Y_s(\omega)-Y_t(\omega)\bigr
\rrvert \leq\cE_s \bigl( \bigl\llvert Y_t
-Y_t(\omega)\bigr\rrvert \bigr) (\omega) + \sup_{\tau\in\cT^s}
\cE_s \bigl((Y_t -X_\tau) \1_{\{
\tau<t\}}
\bigr) (\omega).
\end{equation}
Set $\delta:=\bd[(s,\omega),(t,\omega)]$. Then $\delta\geq t-s$
and we may use Lemma~\ref{leYunifcont} to estimate the first term
in~(\ref{eqproofYContTime}) as
\begin{eqnarray*}
\cE_s \bigl( \bigl\llvert Y_t -Y_t(\omega)
\bigr\rrvert \bigr) (\omega) &=&\sup_{P\in\cP(s,\omega)} E^P
\bigl[\bigl\llvert Y_t^{s,\omega}-Y_t(\omega)\bigr
\rrvert \bigr]
\\
&\leq&\sup_{P\in\cP(s,\omega)} E^P\bigl[\rho_{\cE}
\bigl(\bigl\llVert (\omega\otimes _s B) - \omega\bigr\rrVert
_t\bigr)\bigr]
\\
&\leq&\sup_{P\in\cP(s,\omega)} E^P\bigl[\rho_{\cE}
\bigl(\delta+ \llVert B\rrVert _{t-s}\bigr)\bigr]
\\
&\leq&\sup_{P\in\cP(s,\omega)} E^P\bigl[\rho_{\cE}
\bigl(\delta+ \llVert B\rrVert _{\delta}\bigr)\bigr].
\end{eqnarray*}
To estimate the second term in~(\ref{eqproofYContTime}), let $\tau
\in\cT^s$. As $Y\leq X$ by the definition of $Y$, we deduce
from~(\ref{eqXunifCont}) that
\begin{eqnarray*}
\cE_s \bigl((Y_t -X_\tau) \1_{\{\tau<t\}}
\bigr) (\omega) &\leq& \cE_s \bigl((X_t -X_\tau)
\1_{\{\tau<t\}} \bigr) (\omega)
\\
&=& \sup_{P\in\cP(s,\omega)} E^P \bigl[ \bigl(X_t^{s,\omega}-(X_\tau
)^{s,\omega} \bigr) \1_{\{\tau^{s,\omega}<t\}} \bigr]
\\
&\leq&\sup_{P\in\cP(s,\omega)}
E^P \bigl[\rho_X \bigl(\bd \bigl[(t,\omega
\otimes_s B), (s,\omega\otimes_s B)\bigr] \bigr) \bigr]
\\
&=& \sup_{P\in\cP(s,\omega)} E^P \bigl[
\rho_X \bigl(\llvert t-s\rrvert + \llVert B\rrVert _{t-s}
\bigr) \bigr]
\\
&\leq&\sup_{P\in\cP(s,\omega)} E^P \bigl[\rho_X
\bigl(\delta+ \llVert B\rrVert _{\delta} \bigr) \bigr].
\end{eqnarray*}
Setting $\rho'=\rho_X\vee\rho_{\cE}$, we conclude that
\[
\bigl\llvert Y_s(\omega)-Y_t(\omega)\bigr\rrvert \leq 2
\sup_{P\in\cP(s,\omega)} E^P \bigl[\rho' \bigl(
\delta+ \llVert B\rrVert _\delta \bigr) \bigr] \leq2\rho(\delta),
\]
where $\rho$ is given by Assumption~\ref{asintegrability}. It
remains to set $\rho_Y=2\rho$.
\end{pf}

\begin{remark}\label{rkcadlagAndHitting}
We see from Lemmas~\ref{leYunifcont} and~\ref{leYContTime} that $Y$
is an adapted process with continuous paths. In view of Assumption~\ref
{asXunifCont}, it follows that $X-Y$ is a
c\`adl\`ag adapted process with nonpositive jumps. This implies that
for every $\varepsilon\geq0$, the hitting time $\tau^\varepsilon
=\inf\{t\in[0,T]\dvtx  X_t-Y_t \leq\varepsilon\}$ coincides with the
contact time
$\inf\{t\in[0,T]\dvtx  (X_t\wedge X_{t-}) -Y_t \leq\varepsilon\}$ and,
therefore, is a stopping time (of\vspace*{1pt} the raw filtration $\F$). Moreover,
it implies the pointwise convergence $\tau^\varepsilon\to\tau
^0\equiv\tau^*$ for $\varepsilon\to0$, which we will also find
useful below.
\end{remark}

Lemmas~\ref{leYunifcont} and~\ref{leYContTime} also yield the
following joint continuity property.

\begin{corollary}\label{coYjointlyCont}
There exists a modulus of continuity $\rho_Y$ such that
\[
\bigl\llvert Y_s(\omega)-Y_t\bigl(\omega'
\bigr)\bigr\rrvert \leq\rho_Y\bigl(\bd\bigl[(s,\omega),(t,\omega )
\bigr]\bigr)+ \rho_{\cE}\bigl(\bigl\llVert \omega-\omega'
\bigr\rrVert _T\bigr)
\]
for all $s, t \in[0,T]$ and $\omega,\omega'\in\Omega$. In
particular, if $\theta\dvtx \Omega\to[0,T]$ is any $\llVert  \cdot\rrVert
_T$-continuous function, then
$Y_\theta$ is again continuous.
\end{corollary}

The following submartingale property is a consequence of the dynamic
programming principle of Lemma~\ref{leDPPdetTimes} and ``optional sampling.''

\begin{lemma}\label{lesubmartStop}
Let $s\in[0,T]$ and $\tau\in\cT_s$. Then
%
\begin{equation}
\label{eqsubmartStop} Y_s \leq\cE_s(Y_\tau).
\end{equation}
\end{lemma}

\begin{pf}
By Lemma~\ref{leDPPdetTimes}, we have $Y_s \leq\cE_s(Y_t)$ for any
deterministic time \mbox{$t\in\cT_s$}.

\begin{longlist}[\emph{Step} 1.]
\item[\emph{Step} 1.] We show that~(\ref{eqsubmartStop})
holds when $\tau\in\cT_s$ has finitely many values $t_1<t_2<\cdots<t_n$.

We proceed by induction. If $n=1$, we are in the deterministic case.
Suppose that the result holds for $n-1$ values; in particular, for the
stopping time $\tau\vee t_2\in\cT_{t_1}$. Then, using the tower property,
\begin{eqnarray*}
\cE_s(Y_\tau) &=& \cE_s (Y_{t_1}
\1_{\{\tau=t_1\}} + Y_{\tau\vee t_2}\1_{\{\tau
>t_1\}} )
\\
&=& \cE_s \bigl(Y_{t_1}\1_{\{\tau=t_1\}} +
\cE_{t_1}(Y_{\tau\vee
t_2})\1_{\{\tau>t_1\}} \bigr)
\\
&\geq&\cE_s (Y_{t_1}\1_{\{\tau=t_1\}} + Y_{t_1}
\1_{\{\tau>t_1\}
} )
\\
&\geq& Y_s.
\end{eqnarray*}
\end{longlist}

\begin{longlist}[\emph{Step} 2.]
\item[\emph{Step} 2.]
Let $\tau\in\cT_s$ be arbitrary.
Let $\tau_n=\inf\{t\in D_n\dvtx   t \geq\tau\}$, where $D_n=\{
k2^{-n}T\dvtx   k=0,\ldots,2^n\}$ for $n\geq1$.
Then each $\tau_n$ is a stopping time with finitely many values, and hence
%
\begin{equation}
\label{eqproofSubmartStopStep1} Y_s \leq\cE_s(Y_{\tau_n}),\qquad n
\geq1
\end{equation}
by step~1. In view of $\llvert  \tau_n - \tau\rrvert   \leq2^{-n}T$, Lemma~\ref
{leYContTime} yields that
\[
\llvert Y_{\tau_n}-Y_{\tau}\rrvert \leq\rho_Y
\bigl(\bd\bigl[(\tau,B), (\tau_n,B)\bigr] \bigr) \leq
\rho_Y \bigl( 2^{-n}T + \bigl\llVert B^{\tau}\bigr
\rrVert _{2^{-n}T} \bigr).
\]
In particular,
\begin{eqnarray*}
\bigl(\llvert Y_{\tau_n}- Y_{\tau}\rrvert \bigr)^{s,\omega}
&\leq& \rho_Y \bigl( 2^{-n}T + \bigl\llVert (\omega
\otimes_s B)^{\tau^{s,\omega
}}\bigr\rrVert _{2^{-n}T} \bigr)
\\
&=&
\rho_Y \bigl( 2^{-n}T + \bigl\llVert B^{\tau^{s,\omega}-s}\bigr
\rrVert _{2^{-n}T} \bigr)
\end{eqnarray*}
and thus
%
\begin{eqnarray}\label{eqYestimateDiffStop}
\bigl\llvert \cE_s(Y_{\tau_n}) (\omega) -
\cE_s(Y_{\tau}) (\omega)\bigr\rrvert &\leq&\cE_s
\bigl(\llvert Y_{\tau_n} - Y_{\tau}\rrvert \bigr) (\omega)
\nonumber
\\
&=& \sup_{P\in\cP(s,\omega)} E^P\bigl[\bigl(\llvert
Y_{\tau_n}- Y_{\tau
}\rrvert \bigr)^{s,\omega}\bigr]
\\
&\leq&\sup_{P\in\cP(s,\omega)} E^P\bigl[\rho_Y
\bigl( 2^{-n}T + \bigl\llVert B^{\tau
^{s,\omega}-s}\bigr\rrVert
_{2^{-n}T}\bigr)\bigr].\nonumber
\end{eqnarray}
Note that $\tau^{s,\omega}-s$ is a stopping time as $\tau\in\cT
_s$. Thus, the right-hand side tends to zero as $n\to\infty$, by
Assumption~\ref{asintegrability}.
In view of~(\ref{eqproofSubmartStopStep1}), this completes the proof.\quad\qed
\end{longlist}\noqed
\end{pf}

Next, we discuss the specifics of the discrete situation as introduced
in Remark~\ref{rkoptimalDiscrete}; recall that we use the same
notation $Y$ for the corresponding value function.

\begin{lemma}\label{lerecursion}
Let $\T=\{t_0,t_1,\ldots,t_n\}$, where $n\in\N$ and $t_0<\cdots
<t_n=T$. Then $Y$ is given by
$Y_{t_n}=X_{t_n}$ and
%
\begin{equation}
\label{eqrecursion} Y_{t_i} = X_{t_i} \wedge\cE_{t_i}(Y_{t_{i+1}}),
\qquad i=0,\ldots,n-1.
\end{equation}
Let $\tau^*=\inf\{t\in\T\dvtx   Y_t=X_t\}$, then $Y_{\cdot\wedge\tau
^*}$ is an $\cE$-martingale on $\T$; that is,
%
\begin{equation}
\label{eqdiscreteMart} Y_{t_i\wedge\tau^*} =\cE_{t_i}(Y_{t_{i+1}\wedge\tau^*}),\qquad
i=0,\ldots, n-1,
\end{equation}
and in particular $Y_0= \cE(X_{\tau^*})$.
\end{lemma}

\begin{pf}
Note that $X_T=Y_T$ by the definition of $Y$.
Let $i<n$. From (the obvious discrete version of) Lemma~\ref{leDPPdetTimes},
\[
Y_{t_i} = \inf_{\tau\in\cT^{t_i}(\T)} \cE_{t_i} (
X_\tau\1 _{\{\tau< t_{i+1}\}} + Y_{t_{i+1}} \1_{\{\tau\geq t_{i+1}\}} ).
\]
For any $\tau\in\cT^{t_i}(\T)$, we have either $\tau\equiv t_i$ or
$\tau\geq t_{i+1}$ identically;
hence, the right-hand side equals $X_{t_i} \wedge\cE
_{t_i}(Y_{t_{i+1}})$, which yields~(\ref{eqrecursion}).

We turn to the martingale property. Let $i<n$. On $\{t_i\geq\tau^*\}
$, we have $Y_{t_{i+1}\wedge\tau^*}=Y_{t_i\wedge\tau^*}$ and hence
$Y_{t_i\wedge\tau^*} =\cE_{t_i}(Y_{t_{i+1}\wedge\tau^*})$; whereas
on $\{t_i<\tau^*\}$, we have $Y_{t_i}<X_{t_i}$ and so~(\ref
{eqrecursion}) yields that
\[
Y_{t_i\wedge\tau^*}= Y_{t_i} = \cE_{t_i}(Y_{t_{i+1}}) =
\cE _{t_i}(Y_{t_{i+1}\wedge\tau^*}).
\]
This completes the proof of~(\ref{eqdiscreteMart}), which by the
tower property, also shows that $Y_0=\cE(Y_{T\wedge\tau^*})=\cE
(Y_{\tau^*})=\cE(X_{\tau^*})$.
\end{pf}

We can now prove the optimality of $\tau^*$ by approximating the
continuous problem with suitable discrete ones.

\begin{lemma}\label{leoptimalTau}
Let $\tau^*=\inf\{t\in[0,T]\dvtx   Y_t=X_t\}$. Then $Y_0=\cE(X_{\tau^*})$.
\end{lemma}

\begin{pf}
For $n\geq1$, let $\T_n=D_n=\{k2^{-n}T\dvtx   k=0,\ldots,2^n\}$. Given
$t\in[0,T]$, we denote by $\cT^t_n:=\cT^t(\T_n)$ the corresponding
set of stopping times and by
\[
Y^n_t:=\inf_{\tau\in\cT^t_n}
\cE_t(X_\tau)
\]
the corresponding value function. In view of $\cT^t_n\subseteq\cT
^t$, we have
\[
Y^n\geq Y\qquad\mbox{on }[0,T]\times\Omega.
\]

\begin{longlist}[\emph{Step} 1.]
\item[\emph{Step} 1.]
There exists a modulus of continuity $\rho$
such that
%
\begin{equation}
\label{eqYclose} \bigl\llvert Y^n - Y\bigr\rrvert \leq\rho
\bigl(2^{-n}\bigr)\qquad\mbox{on }\T_n\times\Omega.
\end{equation}

Indeed, let $n\geq1$, $t\in\T_n$ and $\tau\in\cT^t$. Then
$\vartheta:=\inf\{ t\in\T_n\dvtx   t\geq\tau\}$ is in $\cT^t_n$ and
$0\leq\vartheta-\tau\leq2^{-n}T$.
Therefore, Assumption~\ref{asXunifCont} yields that
\begin{eqnarray*}
(X_{\vartheta}- X_{\tau})^{t,\omega} &\leq&\rho_X
\bigl( \bd\bigl[ \bigl(\vartheta^{t,\omega},\omega\otimes_t B
\bigr), \bigl(\tau^{t,\omega},\omega\otimes_t B\bigr)\bigr] \bigr)
\\
&\leq&\rho_X \bigl( 2^{-n}T + \bigl\llVert
B^{\tau^{t,\omega}-t}\bigr\rrVert _{\vartheta
^{t,\omega}-\tau^{t,\omega}} \bigr)
\\
&\leq&\rho_X \bigl( 2^{-n}T + \bigl\llVert
B^{\theta}\bigr\rrVert _{2^{-n}T} \bigr),
\end{eqnarray*}
where $\theta:=\tau^{t,\omega}-t\in\cT$, and hence
\begin{eqnarray*}
\cE_t(X_{\vartheta}) (\omega) - \cE_t(X_{\tau})
(\omega) &\leq&\cE_t(X_{\vartheta}-X_{\tau}) (\omega)
\\
&\leq&\sup_{P\in\cP(t,\omega)} E^P \bigl[\rho_X
\bigl( 2^{-n}T + \bigl\llVert B^{\theta}\bigr\rrVert
_{2^{-n}T} \bigr) \bigr]
\\
&\leq&\rho\bigl(2^{-n}\bigr)
\end{eqnarray*}
for some modulus of continuity $\rho$, by Assumption~\ref
{asintegrability}. As a result,
\[
0 \leq Y^n_t - Y_t = \inf
_{\vartheta\in\cT^t_n} \cE_t(X_\vartheta ) - \inf
_{\tau\in\cT^t} \cE_t(X_\tau) \leq\rho
\bigl(2^{-n}\bigr).
\]
\end{longlist}

\begin{longlist}[\emph{Step} 2.]
\item[\emph{Step} 2.]
Fix $\varepsilon>0$ and define
$\tau^\varepsilon=\inf\{t\in[0,T]\dvtx  X_t-Y_t \leq\varepsilon\}$.
There exists a modulus of continuity $\rho'$ such that for all $n$ satisfying
$\rho(2^{-n})< \varepsilon$,
\[
Y_0 \geq\cE(Y_{\tau^\varepsilon}) + 2\rho\bigl(2^{-n}\bigr) +
\rho'\bigl(2^{-n}\bigr).
\]

Indeed, let $\tau^*_n=\inf\{t\in\T_n\dvtx  Y^n_t=X_t\}$. As $\rho
(2^{-n})< \varepsilon$,~(\ref{eqYclose}) entails that
\[
X-Y^n>0\qquad\mbox{on } \bigl[\!\bigl[0,\tau^{\varepsilon}\bigr[\!\bigr[ \cap (
\T_n \times\Omega);
\]
that is, we have $\tau^{\varepsilon}\leq\tau^*_n$. Define the
stopping time
\[
\tau^{\varepsilon,n}=\inf\bigl\{t\in\T_n\dvtx  t\geq\tau^\varepsilon
\bigr\}.
\]
Recalling that $\tau^*_n$ takes values in $\T_n$, we see that $\tau
^{\varepsilon}\leq\tau^*_n$ even implies that
\[
\tau^{\varepsilon,n}\leq\tau^*_n.
\]
By Lemma~\ref{lerecursion}, the process $Y^n_{\cdot\wedge\tau
^*_n}$ is an $\cE$-martingale on $\T_n$; in particular, using an
optional sampling argument as in step~1 of the proof of Lemma~\ref
{lesubmartStop},
\[
Y^n_0 = \cE\bigl( Y^n_{\tau^{\varepsilon,n}\wedge\tau^*_n}\bigr)
= \cE\bigl( Y^n_{\tau^{\varepsilon,n}}\bigr).
\]
In view of~(\ref{eqYclose}), this implies that
%
\begin{equation}
\label{eqproofOptimalityEpsn} Y_0 \geq \cE( Y_{\tau^{\varepsilon,n}}) - 2\rho
\bigl(2^{-n}\bigr).
\end{equation}
On the other hand, an estimate similar to~(\ref{eqYestimateDiffStop})
and Assumption~\ref{asintegrability} entail that
\[
\bigl\llvert \cE(Y_{\tau^{\varepsilon,n}}) - \cE(Y_{\tau^\varepsilon})\bigr\rrvert \leq
\sup_{P\in\cP} E^P\bigl[\rho_Y\bigl(
2^{-n}T + \bigl\llVert B^{\tau^\varepsilon}\bigr\rrVert _{2^{-n}T}
\bigr)\bigr]\leq\rho'\bigl(2^{-n}\bigr)
\]
for some modulus of continuity $\rho'$. Together with~(\ref
{eqproofOptimalityEpsn}), this yields the claim.
\end{longlist}

\begin{longlist}[\emph{Step} 3.]
\item[\emph{Step} 3.]
Letting $n\to\infty$, step~2
implies that
\[
Y_0\geq\cE(Y_{\tau^\varepsilon}).
\]
%
By Remark~\ref{rkcadlagAndHitting}, we have $\tau^\varepsilon\to
\tau^*$ for $\varepsilon\to0$, and as $Y$ has continuous paths
(Lemma~\ref{leYContTime}), it follows that $Y_{\tau^\varepsilon}\to
Y_{\tau^*}$ pointwise. Thus, (an obvious version of) Fatou's lemma
yields that $Y_0\geq\cE(Y_{\tau^*})$. Recalling the definition of
$\tau^*$, we conclude that
\[
Y_0\geq\cE(Y_{\tau^*})=\cE(X_{\tau^*})\geq\inf
_{\tau\in\cT} \cE(X_\tau) = Y_0.
\]
This completes the proof.\quad\qed
\end{longlist}\noqed
\end{pf}

\begin{remark}\label{rksupermart}
The process $Y_{\cdot\wedge\tau^*}$ is a $P$-supermartingale for any
$P\in\cP$.
\end{remark}
\begin{pf}
Let $0\leq s\leq t\leq T$, where $s\in\T_n$ for some $n$, and let
$\tau\in\cT$ be such that $\tau\leq\tau^\varepsilon$. Going
through step~2 of the preceding proof with the appropriate
modifications then shows that $Y_{s\wedge\tau}\geq\cE_s(Y_\tau)$.
Fix\vspace*{1pt} $P\in\cP$ and note that Assumption~\ref{asanalytic}(ii) implies
$\cE_s(Y_{\tau})\geq E^P[Y_{\tau}| \cF_s]$ $P$-a.s. Choosing $\tau
=t\wedge\tau^\varepsilon$, we obtain that $Y_{s\wedge\tau
^\varepsilon}\geq E^P[Y_{t\wedge\tau^\varepsilon}| \cF_s]$ $P$-a.s.
Now let $\varepsilon\to0$, then Fatou's lemma yields that $Y_{s\wedge
\tau^*}\geq E^P[Y_{t\wedge\tau^*}| \cF_s]$. This shows that
$Y_{\cdot\wedge\tau^*}$ is a $P$-supermartingale on $\bigcup_n\T_n$,
and as $Y$ is continuous, this implies the claim.
\end{pf}

\subsection{Existence of the value}

The aim of this subsection is to show that the upper value function $Y$
coincides with the lower one, denoted by $Z$ below. As mentioned in the
\hyperref[seintro]{Introduction}, there is an obstruction to directly proving the dynamic
programming principle for $Z$ in continuous time; namely, we are unable
to perform measurable selections on the set of stopping times in the
absence of a reference measure. This is related to the measurability
problems that are well known in the literature; see, for example, \cite{FlemingSouganidis89}. In the following lemma, we consider the
discrete setting and prove simultaneously the dynamic programming for
the lower value and that it coincides with the upper value.

\begin{lemma}\label{lerecursionZ}
Let $\T=\{t_0,t_1,\ldots,t_n\}$, where $n\in\N$ and $t_0<\cdots
<t_n=T$, define the lower value function
\[
Z_{t_i}(\omega):= \sup_{P\in\cP(t_i,\omega)} \inf
_{\tau\in\cT
^{t_i}(\T)} E^P\bigl[(X_\tau)^{t_i,\omega}
\bigr],\qquad i=0,\ldots,n,
\]
and recall the upper value function $Y$ introduced in Remark~\ref
{rkoptimalDiscrete}.
For any $j=0,\ldots,n$, we have
%
\begin{equation}
\label{eqZequalsYdiscrete} Z_{t_j}=Y_{t_j}
\end{equation}
and
%
\begin{equation}
\label{eqdppZ} Z_{t_i}(\omega) = \sup_{P\in\cP(t_i,\omega)} \inf
_{\tau\in\cT
^{t_i}(\T)} E^P \bigl[ (X_\tau
\1_{\{\tau< t_{j}\}} + Z_{t_{j}}\1 _{\{\tau\geq t_{j}\}} )^{t_i,\omega} \bigr]
\end{equation}
for all $i=0,\ldots,j$ and $\omega\in\Omega$.
\end{lemma}

\begin{pf}
We proceed by backward induction over $j$. As
$Z_{t_n}=X_{t_n}=Y_{t_n}$, the claim is clear for $j=n$; we show the
passage from $j+1$ to $j$. That is, we assume that for some fixed
$j<n$, we have
%
\begin{equation}
\label{eqZequalsYdiscreteInd} Z_{t_{j+1}}=Y_{t_{j+1}}
\end{equation}
(which, in particular, entails that $Z_{t_{j+1}}$ is $\cF
_{t_{j+1}}$-measurable) and
%
\begin{eqnarray}\label{eqdppZjplus1}
Z_{t_i}(\omega) = \sup_{P\in\cP(t_i,\omega)} \inf
_{\tau\in\cT
^{t_i}(\T)} E^P \bigl[ (X_\tau
\1_{\{\tau< t_{j+1}\}} + Z_{t_{j+1}}\1_{\{\tau\geq t_{j+1}\}} )^{t_i,\omega} \bigr],
\nonumber\\[-8pt]\\[-10pt]
\eqntext{i=0,\ldots,j+1}
\end{eqnarray}
for all $\omega\in\Omega$.
We first note that if $\tau\in\cT^{t_j}(\T)$, then either $\tau
\equiv t_j$ or $\tau>t_j$ identically; therefore,~(\ref
{eqdppZjplus1}) yields that
\begin{eqnarray*}
Z_{t_j}(\omega) &=& \sup_{P\in\cP(t_j,\omega)} \inf
_{\tau\in\cT^{t_j}(\T)} E^P \bigl[ (X_\tau
\1_{\{\tau< t_{j+1}\}} + Z_{t_{j+1}}\1_{\{\tau
\geq t_{j+1}\}} )^{t_j,\omega} \bigr]
\\
&=& X_{t_j}(\omega) \wedge\sup_{P\in\cP(t_j,\omega)} E^P
\bigl[Z_{t_{j+1}}^{t_j,\omega}\bigr]
\\
&=& X_{t_j}(\omega) \wedge\cE_{t_j}(Z_{t_{j+1}}) (
\omega).
\end{eqnarray*}
By the induction assumption~(\ref{eqZequalsYdiscreteInd}) and the
recursion~(\ref{lerecursion}) for $Y$, this shows that
\[
Z_{t_j}= X_{t_j}\wedge\cE_{t_j}(Z_{t_{j+1}})
= X_{t_j}\wedge\cE _{t_j}(Y_{t_{j+1}})=Y_{t_j},
\]
which is~(\ref{eqZequalsYdiscrete}). In particular, $Z_{t_j}$ is $\cF
_{t_j}$-measurable.\vspace*{1pt}

Let us now fix $i\in\{0,\ldots,j\}$ and prove the remaining
claim~(\ref{eqdppZ}). To this end, we first rewrite the latter equation:
substituting the just obtained expression
$Z_{t_j} = X_{t_j} \wedge\cE_{t_j}(Z_{t_{j+1}})$ for $Z_{t_j}$ in the
right-hand side of~(\ref{eqdppZ}), and using~(\ref{eqdppZjplus1})
to substitute $Z_{t_i}$ on the left-hand side of~(\ref{eqdppZ}), we
see that our claim is equivalent to the identity
%
\begin{eqnarray}
\label{eqdppZequiv}
\qquad && \sup_{P\in\cP(t_i,\omega)} \inf_{\tau\in\cT^{t_i}(\T)}
E^P \bigl[ (X_\tau\1_{\{\tau< t_{j+1}\}} + Z_{t_{j+1}}
\1_{\{\tau\geq
t_{j+1}\}} )^{t_i,\omega} \bigr]
\nonumber\\[-8pt]\\[-8pt]\nonumber
&&\qquad = \sup_{P\in\cP(t_i,\omega)} \inf_{\tau\in\cT^{t_i}(\T)} E^P
\bigl[ \bigl(X_\tau\1_{\{\tau< t_{j}\}} + \bigl\{X_{t_j}\wedge
\cE _{t_j}(Z_{t_{j+1}})\bigr\}\1_{\{\tau\geq t_{j}\}}
\bigr)^{t_i,\omega} \bigr].
\end{eqnarray}

We first show the inequality ``$\geq$'' in this equation. To this end,
let $\omega\in\Omega$, $\tau\in\cT^{t_i}(\T)$ and $P\in\cP
(t_i,\omega)$. In view of~(\ref{eqZequalsYdiscreteInd}), Lemma~\ref
{leYunifcont} yields that $Z_{t_{j+1}}$ is continuous and in
particular upper semianalytic. Given $\varepsilon>0$, it then follows
from Assumption~\ref{asanalytic} and an application of the
Jankov--von Neumann selection theorem similar to step~2 of the proof
of~\cite{NutzVanHandel12}, Theorem~2.3, that there exists an \mbox{$\cF
_{t_j-t_i}$-}measurable kernel $\nu\dvtx \Omega\to\mathfrak{P}(\Omega)$
such that
%
\begin{equation}
\label{eqepsOptimal} \qquad\nu(\cdot)\in\cP(t_j,\omega\otimes_{t_i}
\cdot)\quad\mbox{and}\quad E^{\nu(\cdot)} \bigl[Z_{t_{j+1}}^{t_j,\omega\otimes
_{t_i}\cdot}
\bigr] \geq\cE_{t_j}(Z_{t_{j+1}})^{t_i,\omega}(\cdot )-
\varepsilon
\end{equation}
hold $P$-a.s.
Let $\bar{P}$ be the measure defined by
\[
\bar{P}(A)=\iint(\1_A)^{t_j-t_i,\omega'}\bigl(\omega''
\bigr) \nu\bigl(d\omega '';\omega'\bigr)
P\bigl(d\omega'\bigr),\qquad A\in\cF;
\]
then $\bar{P}\in\cP(t_i,\omega)$ by Assumption~\ref
{asanalytic}(iii); moreover, $\bar{P}^{t_j-t_i,\cdot}=\nu(\cdot)$
$P$-a.s. and $\bar{P}=P$ on $\cF_{t_j-t_i}$. In view of~(\ref
{eqepsOptimal}), we have
\begin{eqnarray*}
E^{\bar{P}} \bigl[Z^{t_i,\omega}_{t_{j+1}}| \cF_{t_j-t_i}
\bigr](\cdot) &=& E^{\bar{P}^{t_j-t_i,\cdot}} \bigl[Z^{t_j,\omega\otimes
_{t_i}\cdot}_{t_{j+1}} \bigr]
\\
&=& E^{\nu(\cdot)} \bigl[Z^{t_j,\omega\otimes_{t_i}\cdot
}_{t_{j+1}} \bigr]
\\
&\geq&\cE_{t_j}(Z_{t_{j+1}})^{t_i,\omega}(\cdot)-\varepsilon
\qquad P\mbox{-a.s.}
\end{eqnarray*}
Using this inequality and the tower property of $E^P[ \cdot ]$, we
deduce that
\begin{eqnarray*}
&& E^P \bigl[ \bigl(X_\tau \1_{\{\tau< t_{j}\}} + \bigl
\{X_{t_j}\wedge\cE _{t_j}(Z_{t_{j+1}})\bigr\}
\1_{\{\tau\geq t_{j}\}} \bigr)^{t_i,\omega} \bigr]
\\
&&\qquad \leq E^P \bigl[ \bigl(X_\tau\1_{\{\tau< t_{j}\}}
+X_{t_j}\1_{\{\tau=
t_{j}\}}+ \cE_{t_j}(Z_{t_{j+1}})
\1_{\{\tau\geq t_{j+1}\}} \bigr)^{t_i,\omega} \bigr]
\\
&&\qquad = E^P \bigl[ \bigl(X_\tau\1_{\{\tau< t_{j+1}\}} + \cE
_{t_j}(Z_{t_{j+1}})\1_{\{\tau\geq t_{j+1}\}} \bigr)^{t_i,\omega}
\bigr]
\\
&&\qquad \leq E^{\bar{P}} \bigl[ (X_\tau\1_{\{\tau< t_{j+1}\}} +
Z_{t_{j+1}}\1_{\{\tau\geq t_{j+1}\}} )^{t_i,\omega} \bigr] +\varepsilon.
\end{eqnarray*}
As $\varepsilon>0$, $\tau\in\cT^{t_i}(\T)$ and $P\in\cP
(t_i,\omega)$ were arbitrary, this implies the inequality ``$\geq$''
in~(\ref{eqdppZequiv}).

It remains to show the inequality ``$\leq$'' in~(\ref{eqdppZequiv}).
To this end, let $\omega\in\Omega$, $\tau\in\cT^{t_i}(\T)$,
$P\in\cP(t_i,\omega)$ and define
\[
\btau:=\tau\1_{\{\tau< t_{j}\}} + (t_{j}\1_\Lambda+
t_{j+1}\1 _{\Lambda^c})\1_{\{\tau\geq t_{j}\}},\qquad\Lambda:=\bigl
\{X_{t_{j}} \leq \cE_{t_{j}}(Z_{t_{j+1}})\bigr\}.
\]
Noting that $Z_{t_j}= X_{t_j}\wedge\cE_{t_j}(Z_{t_{j+1}})$ yields
$\Lambda=\{X_{t_j}=Z_{t_j}\}\in\cF_{t_j}$, we see that $\btau\in
\cT_{t_i}(\T)$.
After observing that~(\ref{eqdefE}) and Assumption~\ref
{asanalytic}(ii) imply
\[
\cE_{t_{j}}(Z_{t_{j+1}})^{t_i,\omega} \geq E^P
\bigl[Z_{t_{j+1}}^{t_i,\omega}| \cF_{t_j-t_i}\bigr]\qquad P
\mbox{-a.s.},
\]
we can then use the tower property of $E^P[ \cdot ]$ to obtain that
\begin{eqnarray*}
&& E^P \bigl[ (X_{\btau}\1_{\{\btau< t_{j+1}\}} + Z_{t_{j+1}}
\1_{\{
\btau\geq t_{j+1}\}} )^{t_i,\omega} \bigr]
\\
&&\qquad = E^P \bigl[ \bigl(X_{\tau}\1_{\{\tau< t_{j}\}} +
\{X_{t_j}\1 _{\Lambda} + Z_{t_{j+1}}\1_{\Lambda^c}\}
\1_{\{\tau\geq t_{j}\}} \bigr)^{t_i,\omega} \bigr]
\\
&&\qquad \leq E^P \bigl[ \bigl(X_{\tau}\1_{\{\tau< t_{j}\}} + \bigl
\{X_{t_j}\1 _{\Lambda} + \cE_{t_{j}}(Z_{t_{j+1}})
\1_{\Lambda^c}\bigr\} \1_{\{\tau
\geq t_{j}\}} \bigr)^{t_i,\omega} \bigr]
\\
&&\qquad = E^P \bigl[ \bigl(X_{\tau}\1_{\{\tau< t_{j}\}} + \bigl
\{X_{t_j}\wedge\cE _{t_{j}}(Z_{t_{j+1}})\bigr\}
\1_{\{\tau\geq t_{j}\}} \bigr)^{t_i,\omega
} \bigr].
\end{eqnarray*}
As $\tau\in\cT^{t_i}(\T)$ and $P\in\cP(t_i,\omega)$ were
arbitrary, this implies the desired inequality ``$\leq$'' in~(\ref
{eqdppZequiv}). Here, we have used the fact that, similarly as in
Lemma~\ref{leYaltDef}, the left-hand side of~(\ref{eqdppZequiv})
does not change if we replace $\cT^{t_i}(\T)$ by $\cT_{t_i}(\T)$.
\end{pf}

We can now show the existence of the value for the continuous-time game
by an approximation argument.

\begin{lemma}\label{leexistenceOfValue}
For all $(t,\omega)\in[0,T]\times\Omega$, we have
\[
Z_t(\omega):= \sup_{P\in\cP(t,\omega)} \inf
_{\tau\in\cT^{t}} E^P\bigl[(X_\tau)^{t,\omega}
\bigr] = \inf_{\tau\in\cT^{t}} \sup_{P\in\cP
(t,\omega)}
E^P\bigl[(X_\tau)^{t,\omega}\bigr] \equiv
Y_t(\omega).
\]
\end{lemma}

\begin{pf}
Let $t\in[0,T]$. The inequality $Z_t\leq Y_t$ is immediate from the
ordering of infima and suprema in the definitions; we prove the reverse
inequality.
Given $n\in\N$, we consider $\T_n=\{t\}\cup\{k2^{-n}T\dvtx   k=0,\ldots,2^n\}$ and denote by $Y^n$ and $Z^n$ the corresponding upper and lower
value functions as in Lemmas~\ref{lerecursion} and~\ref{lerecursionZ}, respectively.
As in step~1 of the proof of Lemma~\ref{leoptimalTau}, there exists a
modulus of continuity $\rho$ such that
\[
\bigl\llvert Y^n-Y\bigr\rrvert \leq\rho\bigl(2^{-n}\bigr)
\qquad\mbox{on }\T_n\times\Omega.
\]
Moreover, a similar argument as in the mentioned step shows that
\[
\bigl\llvert Z^n-Z\bigr\rrvert \leq\rho\bigl(2^{-n}\bigr)
\qquad\mbox{on }\T_n\times\Omega.
\]
Since $Y^n_t=Z^n_t$ by Lemma~\ref{lerecursionZ}, we deduce that
$\llvert  Y_t -Z_t\rrvert  \leq2\rho(2^{-n})$, and now the claim follows by letting
$n$ tend to infinity.
\end{pf}

\subsection{Existence of $P^*$}

An important tool in this subsection is an approximation of hitting
times by continuous random times, essentially taken from~\cite{EkrenTouziZhang12stop}.

\begin{lemma}\label{lecontRandomTimes}
Let $\cP$ be weakly compact and $\tau_n=\inf\{t\dvtx   X_t-Y_t\leq
2^{-n}\}$ for $n\geq1$. There exist continuous, $\cF_T$-measurable
functions $\theta_n\dvtx \Omega\to[0,T]$ and \mbox{$\cF_T$-}measurable sets
$\Omega_n\subseteq\Omega$ such that
\[
\sup_{P\in\cP} P\bigl(\Omega_n^c\bigr)
\leq2^{-n}\quad\mbox{and}\quad\tau _{n-1}-2^{-n}\leq
\theta_{n}\leq\tau_{n+1}+2^{-n}\qquad\mbox{on }\Omega_n.
\]
Moreover, $\theta_n\to\tau^*$ $P$-a.s. for all $P\in\cP$.
\end{lemma}

\begin{pf}
In view of Lemmas~\ref{leYunifcont} and~\ref{leYContTime} and
Assumption~\ref{asXunifCont}, the first claim can be argued like
step~1 in the proof of~\cite{EkrenTouziZhang12stop}, Theorem~3.3.
Since $\tau_n\to\tau^*$ by Remark~\ref{rkcadlagAndHitting} and
$\llvert  \tau_n-\theta_n\rrvert  \leq2^{-n}$ on $\Omega_n$, the second claim
follows via the Borel--Cantelli lemma.
\end{pf}

We first establish a measure $P^\sharp$ whose restriction to $\cF
_{\tau^*}$ will be used in the construction of the saddle point.

\begin{lemma}\label{lesaddle1}
Let $\cP$ be weakly compact. Then there exists $P^\sharp\in\cP$
such that $E^{P^\sharp}[X_{\tau^*}]=Y_0$.
\end{lemma}

\begin{pf}
As $X_{\tau^*}=Y_{\tau^*}$, we need to find $P^\sharp\in\cP$ such
that $E^{P^\sharp}[Y_{\tau^*}]\geq Y_0$; the reverse inequality is
clear from Lemma~\ref{leoptimalTau}. For $n\geq1$, let $\tau_n$ and
$\theta_n$ be as in Lemma~\ref{lecontRandomTimes}. In view of
Lemma~\ref{lesubmartStop} and the definition of $\cE(\cdot)$, there
exist $P_n\in\cP$ such that
%
\begin{equation}
\label{eqPn} E^{P_n}[Y_{\tau_n}] \geq\cE(Y_{\tau_n})-2^{-n}
\geq Y_0 - 2^{-n}.
\end{equation}
By passing to a subsequence, we may assume that $P_n\to P^\sharp$
weakly, for some $P^\sharp\in\cP$. Recall that $Y$ is bounded.
As $\theta_n\to\tau^*$ $P^\sharp$-a.s., it follows from
Corollary~\ref{coYjointlyCont} that $E^{P^\sharp}[Y_{\theta_n}]\to
E^{P^\sharp}[Y_{\tau^*}]$. Moreover, the weak convergence $P_m\to
P^\sharp$ and Corollary~\ref{coYjointlyCont} imply that
$\lim_{m}E^{P_m}[Y_{\theta_n}] = E^{P^\sharp}[Y_{\theta_n}]$ for
any fixed $n$. As a result, we have
%
\begin{equation}
\label{eqPstarOne} E^{P^\sharp}[Y_{\tau^*}] = \lim_n
\lim_m E^{P_m}[Y_{\theta_n}].
\end{equation}
On the other hand, for $m\geq n$, we observe that $\tau_m\geq\tau_n$
and, therefore, $E^{P_m}[Y_{\tau_n}]\geq E^{P_m}[Y_{\tau_m}]$ by the
supermartingale property mentioned in Remark~\ref{rksupermart}. Using
also~(\ref{eqPn}), we deduce that
\[
\liminf_n \liminf_m
E^{P_m}[Y_{\tau_n}] \geq\liminf_m
E^{P_m}[Y_{\tau_m}] \geq Y_0.
\]
Combining this with~(\ref{eqPstarOne}), we obtain that
%
\begin{eqnarray}\label{eqPstarTwo}
Y_0-E^{P^\sharp}[Y_{\tau^*}] &\leq&\liminf
_n \liminf_m E^{P_m}[Y_{\tau_n}]
- \lim_n \lim_m E^{P_m}[Y_{\theta_n}]
\nonumber\\[-8pt]\\[-8pt]\nonumber
&\leq&\limsup_n \limsup_m
E^{P_m}\bigl[\llvert Y_{\tau_n}-Y_{\theta_n}\rrvert
\bigr].\nonumber
\end{eqnarray}
It remains to show that the right-hand side vanishes. To this end, we
first note that Lemma~\ref{lecontRandomTimes} yields
\[
\theta_{n-1} - 2^{1-n} \leq\tau_n \leq
\theta_{n+1} + 2^{-n-1}\qquad\mbox{on }\Omega_{n-1}
\cap\Omega_{n+1}.
\]
Of course, we also have $0\leq\tau_n\leq T$. Thus, setting
\[
\psi_n = \sup \bigl\{\llvert Y_t - Y_{\theta_n}
\rrvert\dvtx  t\in\bigl[\theta_{n-1} - 2^{1-n},
\theta_{n+1} + 2^{-n-1}\bigr] \cap[0,T] \bigr\},
\]
we have
\begin{eqnarray*}
E^{P_m}\bigl[\llvert Y_{\tau_n}-Y_{\theta_n}\rrvert \bigr]
&\leq& E^{P_m}[\psi_n]+ 4\llVert Y\rrVert _{\infty}
P_m\bigl(\Omega_{n-1}^c\cup\Omega
_{n+1}^c\bigr)
\\
&\leq& E^{P_m}[\psi_n]+ 2^{4-n}\llVert Y\rrVert
_{\infty}.
\end{eqnarray*}
Moreover, $\psi_n$ is uniformly bounded, and the continuity of $\theta
_k$ and Corollary~\ref{coYjointlyCont} yield that $\psi_n$ is continuous.
Therefore, $E^{P_m}[\psi_n] \to E^{P^\sharp}[\psi_n]$ for each $n$,
and we conclude that
\begin{eqnarray*}
\limsup_n \limsup_m E^{P_m}
\bigl[\llvert Y_{\tau_n}-Y_{\theta_n}\rrvert \bigr] &\leq&\limsup
_n \limsup_m E^{P_m}[
\psi_n]
\\
&\leq&\limsup_n E^{P^\sharp}[\psi_n] =0,
\end{eqnarray*}
where the last step used dominated convergence under $P^\sharp$ and
the fact that $\psi_n\to0$ $P^\sharp$-a.s. due to $\theta_n\to
\tau^*$ $P^\sharp$-a.s. In view of~(\ref{eqPstarTwo}), this
completes the proof.
\end{pf}

The measure $P^\sharp$ already satisfies
\[
\inf_{\tau\in\cT, \tau\leq\tau^* }E^{P^\sharp}[X_{\tau}] =
E^{P^\sharp}[X_{\tau^*}].
\]
(This follows using Remark~\ref{rksupermart}; cf. the proof of
Lemma~\ref{lesaddle2} below.) In order to obtain a saddle point, we
need to find an extension of $P^\sharp| _{\cF_\tau^*}$ to $\cF$
under which ``after $\tau^*$, immediate stopping is optimal.'' As a
preparation, we first note the following semicontinuity property.

\begin{lemma}\label{levalueUsc}
The function $P\mapsto\inf_{\tau\in\cT^t} E^P[(X_\tau)^{t,\omega
}]$ is upper semicontinuous on $\fPO$, for all $(t,\omega)\in
[0,T]\times\Omega$.
\end{lemma}

\begin{pf}
We state the proof for $t=0$; the general case is proved similarly. Let
$P_n\to P$ in $\fPO$; we need to show that
\[
\limsup_{n\to\infty} \inf_{\tau\in\cT}
E^{P_n}[X_\tau] \leq \inf_{\tau\in\cT}
E^{P}[X_\tau].
\]
To this end, it suffices to show that given $\varepsilon>0$ and $\tau
\in\cT$, there exists $\tau'\in\cT$ such that
%
\begin{equation}
\label{eqproofValueUsc} \limsup_{n\to\infty} E^{P_n}[X_{\tau'}]
\leq E^{P}[X_\tau] + \varepsilon.
\end{equation}
Moreover, by an approximation from the right, we may suppose that $\tau
=\sum_{i=1}^N t_i\1_{A_i}$ for some $N\in\N$, $t_i\in[0,T]$ and
$A_i\in\cF_{t_i}$. Given $\delta>0$, we can find for each $1\leq
i\leq N$ a $P$-continuity set $D_i\in\cF_{t_i}$ [i.e., $P(\partial
D_i)=0$] satisfying $P(A_i \Delta D_i)<\delta$. Note that
$D_0:=(D_1\cup\cdots\cup D_N)^c$ is then also a $P$-continuity set.
Define $t_0:=T$ and
\[
\tau':= \sum_{i=0}^N
t_i\1_{D_i};
\]
then $\tau'\in\cT$ and $P\{\tau\neq\tau'\}<N\delta$. As $X$ is
bounded, it follows that $E^{P}[\llvert  X_\tau-X_{\tau'}\rrvert  ]<\varepsilon$
for $\delta>0$ chosen small enough, while
\[
E^{P_n}[X_{\tau'}] \to E^{P}[X_{\tau'}]
\]
since $X_{\tau'}=\sum_{i=0}^N X_{t_i}\1_{D_i}$, each $X_{t_i}$ is
bounded and continuous, and each $D_i$ is a $P$-continuity set. This
implies~(\ref{eqproofValueUsc}).
\end{pf}

We can now construct the kernel that will be used to extend $P^\sharp$.

\begin{lemma}\label{leworstMeasure}
Let $\cP(t,\omega)$ be weakly compact for all $(t,\omega)\in
[0,T]\times\Omega$ and let $\theta\in\cT$. There exists an $\cF
_{\theta}^*$-measurable kernel $\hat{P}_\theta\dvtx  \Omega\to\fPO$
such that $\hat{P}_\theta(\omega)\in\cP(\theta,\omega)$ and
\[
\inf_{\tau\in\cT^{\theta(\omega)}} E^{\hat{P}_\theta(\omega
)}\bigl[(X_\tau)^{\theta,\omega}
\bigr] = \sup_{P\in\cP(\theta,\omega
)}\inf_{\tau\in\cT^{\theta(\omega)}}
E^P\bigl[(X_\tau)^{\theta,\omega}\bigr]
\]
for all $\omega\in\Omega$.
\end{lemma}

\begin{pf}
For brevity, let us define
\[
V(t,\omega,P):=\inf_{\tau\in\cT^t} E^{P}
\bigl[(X_\tau)^{t,\omega}\bigr].
\]
We first fix $P\in\fPO$ and note that $(t,\omega)\mapsto V(t,\omega,P)$ is Borel. To see this, we first observe that
\[
V(t,\omega,P)=\inf_{\tau\in\cT} E^{P}
\bigl[X_{\tau(\cdot)\vee
t}(\omega\otimes\cdot)\bigr]
\]
by the argument of Lemma~\ref{leYaltDef}. Moreover, let $\cT
'\subseteq\cT$ be a countable set such that for each $\tau\in\cT$
there exist $\tau_n\in\cT'$ satisfying
$\tau_n\downarrow\tau$ $P$-a.s.; for instance, $\cT'$ can be chosen
to consists of stopping times of the form $\sum_{i=1}^N t_i\1_{A_i}$,
where each $t_i$ is dyadic and $A_i$ belongs to a countable collection
generating $\cF_{t_i}$. Then we have
\[
V(t,\omega,P)=\inf_{\tau\in\cT^t} E^{P}
\bigl[(X_\tau)^{t,\omega}\bigr] = \inf_{\tau\in\cT'}
E^{P}\bigl[X_{\tau(\cdot)\vee t}(\omega\otimes \cdot)\bigr]
\]
by dominated convergence and as $(t,\omega)\mapsto E^{P}[X_{\tau
(\cdot)\vee t}(\omega\otimes\cdot)]$ is Borel for every $\tau$ by
Fubini's theorem, it follows that
$(t,\omega)\mapsto V(t,\omega,P)$ is Borel.

On the other hand, we know from Lemma~\ref{levalueUsc} that $P\mapsto
V(t,\omega,P)$ is upper semicontinuous. Together, it follows that $V$
is Borel as a function on $[0,T]\times\Omega\times\fPO$; in
particular, $(\omega,P)\mapsto V_\theta(\omega,P):=V(\theta(\omega
),\omega,P)$ is again Borel (recall that $\cT$ consists of $\F
$-stopping times).

For each $(t,\omega)\in[0,T]\times\Omega$, it follows by
compactness and Lemma~\ref{levalueUsc} that there exists
$\hat{P}\in\cP(t,\omega)$ such that
\[
V(t,\omega,\hat{P}) = \sup_{P\in\cP(t,\omega)} V(t,\omega,P);
\]
in particular, $P\mapsto V_\theta(\omega,P)$ admits a maximizer for
each $\omega\in\Omega$. As the graph of $\cP(\theta,\cdot)$ is
analytic, the Jankov--von Neumann theorem in the form of \cite{BertsekasShreve78}, Proposition~7.50(b), page~184, shows that a
maximizer can be chosen in a universally measurable way, which yields
the claim.
\end{pf}

Finally, we can prove the remaining result of Theorem~\ref{thoptimal}.

\begin{lemma}\label{lesaddle2}
Let $\cP(t,\omega)$ be weakly compact for all $(t,\omega)\in
[0,T]\times\Omega$, let $P^\sharp$ be as in Lemma~\ref{lesaddle1}
and let $\hat{P}_{\tau^*}$ be as in
Lemma~\ref{leworstMeasure}. Then the measure defined by
\[
P^*(A)=\iint(\1_A)^{\tau^*,\omega}\bigl(\omega'\bigr)
\hat{P}_{\tau
^*}\bigl(d\omega';\omega\bigr)
P^\sharp(d\omega),\qquad A\in\cF
\]
is an element of $\cP$ and satisfies
\[
\inf_{\tau\in\cT}E^{P^*}[X_{\tau}] =
E^{P^*}[X_{\tau^*}].
\]
\end{lemma}

\begin{pf}
We set $\hat{P}:=\hat{P}_{\tau^*}$. After replacing $\hat{P}$ with
a Borel kernel $\nu$ such that $\nu=\hat{P}$ $\P^\sharp$-a.s., it
follows from Assumption~\ref{asanalytic}(iii) that $P^*\in\cP$.
Let $\tau\in\cT$; then the definition of $\hat{P}$ and Lemma~\ref
{leexistenceOfValue} yield
\begin{eqnarray*}
E^{\hat{P}(\omega)}\bigl[(X_{\tau\vee\tau^*})^{\tau^*,\omega}\bigr] &\geq&\inf
_{\theta\in\cT^{\tau^*(\omega)}} E^{\hat{P}(\omega
)}\bigl[(X_\theta)^{\tau^*,\omega}
\bigr]
\\
&=& \sup_{P\in\cP(\tau^*,\omega)} \inf_{\theta\in\cT^{\tau
^*(\omega)}} E^{P}
\bigl[(X_\theta)^{\tau^*,\omega}\bigr]
\\
&=& \inf_{\theta\in\cT^{\tau^*(\omega)}} \sup_{P\in\cP(\tau
^*,\omega)} E^{P}
\bigl[(X_\theta)^{\tau^*,\omega}\bigr]
\\
&=& Y_{\tau^*}(\omega)
\end{eqnarray*}
for all $\omega\in\Omega$. This means that $E^{P^*}[X_{\tau\vee
\tau^*}|  \cF_{\tau^*}] \geq Y_{\tau^*}$ $P^*$-a.s., and thus
\begin{eqnarray*}
E^{P^*}[X_\tau| \cF_{\tau^*}] &=& E^{P^*}[X_{\tau\wedge\tau^*}|
\cF_{\tau^*}] \1_{\{\tau< \tau
^*\}} + E^{P^*}[X_{\tau\vee\tau^*}|
\cF_{\tau^*}] \1_{\{\tau\geq
\tau^*\}}
\\
&\geq& X_{\tau\wedge\tau^*} \1_{\{\tau< \tau^*\}} + Y_{\tau^*} \1
_{\{\tau\geq\tau^*\}}\qquad P^*\mbox{-a.s.}
\end{eqnarray*}
By Remark~\ref{rksupermart}, $Y_{\cdot\wedge\tau^*}$ is a
$P^\sharp$-supermartingale, but as $Y_0=E^{P^\sharp}[Y_{\tau^*}]$ by
Lemma~\ref{lesaddle1}, $Y_{\cdot\wedge\tau^*}$ is even a $P^\sharp
$-martingale, and hence a $P^*$-martingale. Using also that $X\geq Y$,
we conclude that
\[
E^{P^*}[X_\tau| \cF_{\tau\wedge\tau^*}] \geq Y_{\tau\wedge\tau^*}
\1_{\{\tau< \tau^*\}} + E^{P^*}[Y_{\tau
^*}| \cF_{\tau\wedge\tau^*}]
\1_{\{\tau\geq\tau^*\}} = Y_{\tau\wedge\tau^*}\qquad P^*\mbox{-a.s.}
\]
and thus
\[
E^{P^*}[X_\tau]\geq E^{P^*}[Y_{\tau\wedge\tau^*}]=E^{P^*}[Y_{\tau
^*}]=E^{P^*}[X_{\tau^*}].
\]
Since $\tau\in\cT$ was arbitrary, this proves the claim.
\end{pf}

\section{Application to American options}\label{seamericanOption}

In this section, we apply our main result to the pricing of American
options under volatility uncertainty.
To this end, we interpret $B$ as the stock price process and assume that
$\cP$ consists of local martingale measures, each of which is seen as
a possible scenario for the volatility. More precisely, following~\cite{SonerTouziZhang2010dual}, we assume that $\cP$ is a subset of $\cP
_S$, the set of all local martingale laws of the form
\[
P^\alpha= P_0 \circ \biggl(\int_0^\cdot
\alpha^{1/2}_u \,dB_u \biggr)^{-1},
\]
where $P_0$ is the Wiener measure and $\alpha$ ranges over all locally
square integrable, progressively measurable processes with values in
$\S_{++}$. We remark that if $\cP$ is not already a subset of $\cP
_S$, then we may replace $\{\cP(s,\omega)\}$ by $\{\cP(s,\omega
)\cap\cP_S\}$ without invalidating Assumption~\ref{asanalytic};
cf. \cite{NeufeldNutz12}, Corollary~2.5.

Let $\G= (\cG_t)_{0\leq t \leq T}$ be the filtration defined by $\cG
_t= \cF^{*}_{t}\vee\mathcal{N}^{\mathcal{P}}$,
where $\cF^{*}_{t}$ is the universal completion of $\cF_t$ and
$\mathcal{N}^{\mathcal{P}}$ is the collection of all sets which are
$(\cF_T,P)$-null for all $P\in\cP$. Let $H$ be an $\R^d$-valued,
$\G$-predictable process such that $\int_0^T H_u^\top \,d\langle B
\rangle_u H_u <\infty$ $P$-a.s. for all $P\in\cP$. Then $H$ is
called an admissible trading strategy if the $P$-integral
$\int H \,dB$ is a $P$-supermartingale, for all $P\in\cP$, and we
denote by $\cH$ the set of all admissible trading strategies.

If $X$ is an American-style option where the buyer chooses the exercise
time, then the buyer's price (or subhedging price) is given by
\begin{eqnarray*}
x_*(X)&:=&\sup \biggl\{x\in\R\dvtx \mbox{there exist }\tau\in\cT\mbox{ and }H\in\cH
\mbox{ such that }
\\
&&\hspace*{36pt} X_\tau+ \int_0^\tau H_u
\,dB_u \geq x\ P\mbox{-a.s. for all }P\in\cP \biggr\}.
\end{eqnarray*}
This is the supremum of all prices $x$ such that, by using a suitable
choice of hedging strategy and exercise time, the buyer will incur no
loss, no matter which scenario $P$ occurs.
On the other hand, if $X$ is a short position in an American option, so
that the seller chooses the exercise time, then the corresponding
sellers's price is given by
\begin{eqnarray*}
x^*(X)&:=&\inf \biggl\{x\in\R\dvtx \mbox{there exist }\tau\in\cT\mbox{ and }H\in\cH
\mbox{ such that }
\\
&&\hspace*{33pt} x + \int_0^\tau H_u
\,dB_u \geq X_\tau\ P\mbox{-a.s. for all }P\in
\cP \biggr\}.
\end{eqnarray*}
Clearly, $x_*(X)=-x^*(-X)$, so it suffices to study one of these cases.
We state the result for the seller's price $x^*$ (because it matches
the sign convention for nonlinear expectations).

\begin{theorem}\label{thduality}
Let Assumptions~\ref{asXunifCont},~\ref{asEunifCont} and~\ref
{asintegrability} hold. Then
\[
x^*(X)=\inf_{\tau\in\cT} \cE(X_{\tau}) =
\cE(X_{\tau^*})\qquad\mbox{for }\tau^*=\inf\bigl\{t\in[0,T]\dvtx
Y_t=X_t\bigr\},
\]
and there exists $H\in\cH$ such that
$x^*(X) + \int_0^\tau H_u \,dB_u \geq X_{\tau^*}$ $P$-a.s. for all
$P\in\cP$; in particular, the infimum defining $x^*(X)$ is attained.
\end{theorem}

\begin{pf}
We set $x^*=x^*(X)$ and $y^*=\inf_{\tau\in\cT} \cE(X_{\tau})$.
Let $x>x^*$, then the definition of $x^*$ yields $\tau\in\cT$ and
$H\in\cH$ such that
\[
x + \int_0^\tau H_u
\,dB_u\geq X_{\tau}\qquad P\mbox{-a.s. for all }P\in
\cP.
\]
As $H$ is admissible, this implies that $x\geq E^P[X_{\tau}]$ for all
$P\in\cP$, and thus $x\geq\cE(X_{\tau})$. In particular, $x\geq
\inf_{\tau\in\cT} \cE(X_\tau)=y^*$. As $x>x^*$ was arbitrary,
this shows that $x^*\geq y^*$.

Conversely, we have $y^*=\cE(X_{\tau^*})$ by Theorem~\ref
{thoptimal}. Moreover, as $X_{\tau^*}$ is Borel-measurable and
bounded, the (European) superhedging result stated in \cite{NeufeldNutz12}, Theorem~2.3, yields $H\in\cH$ such that
\[
\cE(X_{\tau^*}) + \int_0^\tau
H_u \,dB_u\geq X_{\tau^*}\qquad P\mbox{-a.s. for all }P\in\cP.
\]
Thus, the definition of $x^*$ implies that $x^*\leq\cE(X_{\tau^*})=y^*$.
\end{pf}

\begin{remark}
In view of Remark~\ref{rkoptimalDiscrete}, we can show a similar
result for Bermudan options, that is, options where the exercise time
can be chosen from a given set $\T=\{t_0,\ldots,t_n\}$.
\end{remark}\vspace*{-9pt}

\begin{appendix}\label{app}
\section*{Appendix: Proof of Lemma~\texorpdfstring{\lowercase{\protect\ref{lemasssss}}}{3.10}}

In this section, we complete Example~\ref{excontrolledSDE} by showing
that Assumption~\ref{asEunifCont} is satisfied. We use the setting
and notation introduced in that example.

\begin{pf*}{Proof of Lemma~\ref{lemasssss}}
Let $t\in[0,T]$, $\tau\in\cT^t$ and $\bar\omega\in\Omega$.
Using a discretization of stochastic integrals
as in~\cite{Karandikar95} and the fact that the paths of $B$ are
continuous, we can define $\F$-progressively measurable processes
$A^n$ such that
\[
A:= \limsup_{n\to\infty} A^{n}
\]
coincides $P$-a.s. with the usual quadratic variation process of $B$
under $P$, for any semimartingale law $P$. Let us also define (with
$\infty-\infty:=-\infty$, say)
%
\begin{equation}
\label{eqahat} \qquad a_{u}:= \limsup_{n\to\infty} n
(A_{u+1/n} - A_{u} )\quad\mbox{and}\quad\hat{a}_{u}:= a_{u} \1_{\{a_{u}\in\S_{++}\}} + 1\1 _{\{a_{u}\notin\S_{++}\}};
\end{equation}
then $\hat{a}$ is $\F_{+}$-progressively measurable and coincides
$dt\times P$-a.s. with the squared volatility of
$B$ under $P$, for any $P\in\bigcup_{\omega}\cP(t,\omega)$. Given
$\omega\in\Omega$ and recalling that $\sigma$ admits the inverse
$\sigma^{\inv}$ in its third argument, we may then define the
$U$-valued process
\[
\hat{\nu}^{\omega}_{u}:= \sigma^{\inv}\bigl(u+t,
\omega\otimes_{t} \cdot, \hat{a}_{u}^{1/2}\bigr).
\]
Let $\Gamma^{t,\omega,\nu}$ denote the solution of the SDE with
parameters $(t,\omega)$ and control $\nu\in\cU$. For any $\nu\in
\cU$, we have by construction that
\[
\hat{a}\bigl(\Gamma^{t,\omega,\nu}\bigr) = \sigma^{2} \bigl(\cdot+t,
\omega \otimes_{t} \Gamma^{t,\omega,\nu}, \nu\bigr) \qquad
P_{0}\mbox{-a.s.}
\]
and thus
%
\begin{equation}
\label{equnivControlProp} \hat{\nu}^{\omega}\bigl(\Gamma^{t,\omega,\nu}\bigr) = \nu
\qquad P_{0}\mbox{-a.s.}
\end{equation}
We emphasize that these identities indeed hold up to
$P_{0}$-evanescence (rather than just $dt\times P_{0}$-a.s.) because
$\sigma^{2} (\cdot+t,\omega\otimes_{t} \Gamma^{t,\omega,\nu})$
is right-continuous $P_{0}$-a.s. and the ``derivative'' in~(\ref
{eqahat}) is taken from the right. In particular, (\ref
{equnivControlProp}) implies that $\hat{\nu}^{\omega}$ has c\`adl\`ag paths $P(t,\omega,\nu)$-a.s. For later use, we also note that
$(\omega,\omega')\mapsto\hat{\nu}^{\omega}_{u}(\omega')$ is $\cF
_{t}\otimes\cF$-measurable.

Given two paths $\omega,\bar\omega\in\Omega$, let us now consider
the equation
%
\begin{equation}
\label{eqtransformationSDE} \zeta= \int_0^{\cdot} \sigma
\bigl(u+t,\bar\omega\otimes_t \zeta, \hat{\nu}^{\omega}_{u}
\bigr) \sigma\bigl(u+t,\omega\otimes_t B,\hat{\nu }^{\omega}_{u}
\bigr)^{-1} \,dB_u.
\end{equation}
Under $P(t,\omega,\nu)$, there exists an almost-surely unique strong
solution $\zeta^{t,\omega,\nu}$ and it follows via~(\ref
{equnivControlProp}) that $\zeta^{t,\omega,\nu}(\Gamma^{t,\omega,\nu}) = \Gamma^{t,\bar\omega,\nu}$ $P_{0}$-a.s. However, we need
to define the solution universally, without reference to $\nu$. To
this end, we again use a discretization as in~\cite{Karandikar95} to
define approximate solutions $\zeta^n$ (which are \mbox{$\F
_{+}$-}progressively measurable and merely c\`adl\`ag, whence the need
to have $\sigma$ defined on $\mathbbm{D}$) and set $\zeta'^{\omega
}:= \limsup_{n\to\infty} \zeta^{n}$. Since the integrand in~(\ref
{eqtransformationSDE}) is $P(t,\omega,\nu)$-a.s. c\`adl\`ag, we
have that $\zeta'^{\omega}$ coincides with $\zeta^{t,\omega,\nu}$
$P(t,\omega,\nu)$-a.s.; cf.~\cite{Karandikar95}. In particular,
$\zeta'^{\omega}$ is continuous $P(t,\omega,\nu)$-a.s., so that
\[
\zeta^{\omega}_{u}:= \limsup_{q\in\Q,  q\uparrow u} \zeta
'^{\omega}_{q}
\]
still coincides with $\zeta^{t,\omega,\nu}$ $P(t,\omega,\nu
)$-a.s., while in addition being $\F$-progressively measurable.
Moreover, $(\omega,\omega')\mapsto\zeta^{\omega}(\omega')$ is
$\cF_{t}\otimes\cF$-measurable by construction. While we now have
\[
\zeta^{\omega}\bigl(\Gamma^{t,\omega,\nu}\bigr) = \Gamma^{t,\bar\omega,\nu}
\qquad P_{0}\mbox{-a.s.}
\]
simultaneously for all $\nu\in\cU$, as desired, we still have to
elaborate on the definition of $\tau_{\omega}$. Indeed, we cannot
ensure that all paths of $\zeta^{\omega}$ are continuous, so that the
right-hand side of~(\ref{eqdefStoppingTransformSketch}) is not well
defined. (We cannot simply set the irregular paths to zero as in~\cite{Karandikar95}, for then the resulting process would not be $\F
$-adapted and so $\tau_{\omega}$ would not be an $\F$-stopping time
as required.)

To simplify the notation, let $\tilde{\zeta}$ be the process defined by
\[
\tilde{\zeta}^{\omega}\bigl(\omega'\bigr):= 0
\otimes_t \zeta^{\omega
}\bigl(\omega'_{\cdot+t}-
\omega'_t\bigr),\qquad\omega'\in\Omega.
\]
Given $r\in[0,T]$, let $\llVert    \cdot  \rrVert  _{1/3,r}$ be the $1/3$-H\"
older norm for functions considered on $[0,r]\cap\Q$, and note that
its computation involves only a countable supremum. Thus,
\[
C^{\omega}_{r}:= \bigl\{ \omega'\in\Omega\dvtx  \bigl
\llVert \tilde{\zeta }^{\omega}\bigl(\omega'\bigr) \bigr
\rrVert _{1/3,r}<\infty \bigr\} \in\cF_{r}.
\]
Moreover, $\tilde{\zeta}^{\omega}| _{[0,r]}$ is continuous on
$C^{\omega}_{r}$, and a standard result for the path regularity of
martingales shows that $C^{\omega}_r$ has full $P(t,\omega,\nu
)$-measure for any $\nu\in\cU$. Consider
\[
D^{\omega}_{r}:= \bigl\{ \omega'\in
C^{\omega}_{r}\dvtx  \tau \bigl(
\tilde{\zeta}_{\cdot\wedge r}\bigl(\omega'\bigr)\bigr) \leq r \bigr
\} \in\cF_{r}.
\]
By Galmarino's test, we have that $\tau(\tilde{\zeta}_{\cdot\wedge
r})=\tau(\tilde{\zeta}_{\cdot\wedge r'})$ on $D^{\omega}_{r} \cap
D^{\omega}_{r'}$ for any $r,r'\in[0,T]$. Thus,
\[
\tau_{\omega}\bigl(\omega'\bigr):= %
\cases{ \tau
\bigl(\tilde{\zeta}_{\cdot\wedge r}\bigl(\omega'\bigr)\bigr), &\quad
if $\omega'\in D^{\omega}_{r}$, $r\in[0,T]$;
\vspace*{3pt}\cr
T, &\quad if $\displaystyle\omega'\in\biggl(\bigcup_{r}D^{\omega}_{r}
\biggr)^{c}$} %
\]
is well defined. To see that the Borel-measurable function $\tau
_{\omega}$ is an $\F$-stopping time, we observe that $\{\tau_{\omega
}=T\}\in\cF_{T}$ and, for $u<T$,
\[
\{\tau_{\omega}=u\} = \bigl\{ \omega'\in
C^{\omega}_{u}\dvtx  \tau \bigl(\tilde{\zeta}_{\cdot\wedge u}\bigl(
\omega'\bigr)\bigr)=u \bigr\} \in\cF_{u},
\]
due to the fact that $\tilde{\zeta}$ is $\F$-adapted. In fact, we
have $\tau_{\omega}\in\cT^{t}$ by the definition of $\tilde{\zeta
}$ and the condition that $\tau\geq t$. Moreover, $(\omega,\omega
')\mapsto\tau_{\omega}(\omega')$ is $\cF_{t}\otimes\cF
$-measurable by construction.
Since $C^{\omega}_r$ has full $P(t,\omega,\nu)$-measure for any $\nu
\in\cU$, we also have
\[
\tau_\omega\bigl(0\otimes_t \Gamma^{t,\omega,\nu}\bigr) =
\tau\bigl(0\otimes_t \Gamma^{t,\bar\omega,\nu}\bigr)\qquad P_0
\mbox{-a.s.}
\]
for all $\nu\in\cU$, and we deduce as in~(\ref
{eqfinalEstimateSketch}) that
\begin{eqnarray*}
&&\bigl\llvert \cE_t(X_\tau) (\bar\omega) -
\cE_t(X_{\tau_\omega}) (\omega )\bigr\rrvert
\\
&&\qquad \leq \sup_{\nu\in\cU} E^{P_0} \bigl[ \bigl\llvert
X_{\tau(0\otimes_t
\Gamma^{t,\bar\omega,\nu})} \bigl(\bar\omega\otimes_t \Gamma^{t,\bar
\omega,\nu}
\bigr) - X_{\tau(0\otimes_t \Gamma^{t,\bar\omega,\nu})} \bigl(\omega\otimes_t \Gamma^{t,\omega,\nu}
\bigr) \bigr\rrvert \bigr]
\\
&&\qquad \leq\rho_X\bigl(C\llVert \bar\omega-\omega \rrVert _t
\bigr)
\end{eqnarray*}
as desired.
\end{pf*}
\end{appendix}

\section*{Acknowledgments}
The authors thank Ioannis Karatzas and the anonymous referees for their
comments.



\printaddresses
\end{document}